\documentclass[12pt]{article}%
\usepackage{amsmath,amssymb,amsthm,amsfonts}
\usepackage{wasysym}
\usepackage{graphicx}
\usepackage[dvipsnames]{xcolor}
\usepackage{stackengine}

\usepackage[colorlinks]{hyperref}
\usepackage{tikz}
\usepackage[export]{adjustbox}

\usepackage{geometry}
\usepackage{appendix}

\renewcommand{\epsilon}{\varepsilon}

\definecolor{red}{rgb}{0.8500, 0.3250, 0.0980}
\definecolor{green}{rgb}{0.4660, 0.6740, 0.1880}
\definecolor{yellow}{rgb}{0.9290, 0.6940, 0.1250}
\definecolor{blue}{rgb}{0, 0.4470, 0.7410}

\begin{document}

\title{Walking droplets as a damped-driven system}

\author{Aminur Rahman \thanks{Corresponding  \url{arahman2@uw.edu}}
\thanks{Department of Applied Mathematics, University of Washington}, J. Nathan Kutz\footnotemark[2]}

\date{}

\maketitle

\begin{abstract}
We consider the dynamics of a droplet on a vibrating fluid bath.  This  hydrodynamic quantum analog system is shown to elicit the canonical behavior of damped-driven systems, including a period doubling route to chaos.  By approximating the system as a compositional map between the gain and loss dynamics, the underlying nonlinear dynamics can be shown to be driven by energy balances in the systems.  The gain-loss iterative mapping is similar to a normal form encoding for the pattern forming instabilities generated in such spatially-extended system.  Similar to mode-locked lasers and rotating detonation engines, the underlying bifurcations persist for general forms of the loss and gain, both of which admit explicit representations in our approximation.  Moreover, the resulting geometrical description of the particle-wave interaction completely characterizes the instabilities observed in experiments.
\end{abstract}

\bigskip
\bigskip

\section{Introduction}
\label{Sec: Intro}

Over the past decade, walking droplet systems have been studied as hydrodynamic analogs of quantum behavior.  First observed in the seminal works of Couder, Fort, \textit{et al.} \cite{CPFB05, CouderFort06}, such droplet systems exhibit a diversity of dynamical behaviors as the droplet interacts with the vibrating fluid bath. Provided the forcing on the bath is large enough, but below the Faraday wave threshold, the droplet creates persistent wave structures on the bath with each impact.  If the forcing is increased further still, beyond the Faraday wave threshold, standing waves are present even in the absence of a droplet.  These so-called \emph{pilot waves} interact with the droplet and propel the droplet in the horizontal direction, thereby creating a wave-particle coupled system. Much of the known fluid mechanics and quantum analogs have been documented in review articles of Bush \textit{et al.} \cite{Bush10, Bush15a, Bush15b, BushOza20_ROPP, TCB18}. In addition to quantum analogs, the walking droplet system is an example of a damped-driven system for which there exist canonical underlying bifurcations. Indeed, the droplet system can be mathematically quantified by a discrete map accounting for the energy balance in the system, with the driving (gain) being the vibration of the bath and the damping (loss) being modeled by the wave-particle contact dynamics as detailed by Mola\v{c}ek and Bush \cite{MolBush13a, MolBush13b}.  These two gain-loss mechanisms are sufficient to explain the observed bifurcations, or period-doubling route to chaos, of the hydrodynamic analog system.

Much of the hydrodynamic analogs have focused on recreating quantum-scale behavior.  Perhaps some of the most compelling quantum-like observations in walking droplet systems occur on circular corrals \cite{HMFCB13, OHRB14, OWHRB14, Gilet16, CSB18}, elliptical corrals \cite{SCB18}, and with tunneling \cite{EFMC09, NMB17}.  In the hydrodynamic circular corral experiments of Harris \textit{et al.} \cite{HMFCB13} and Cristea-Platon \textit{et al.} \cite{CSB18}, they observe long-time statistics similar to that of the quantum corral experiment \cite{CLE1993}.  In the hydrodynamic elliptical corral studied by S\'{a}enz \textit{et al.} \cite{SCB18} they observe effects similar to that of the quantum mirage \cite{QuantumMirage}.  While the hydrodynamic circular corral phenomena has been studied analytically, the hydrodynamic elliptical corral has yet to be modeled.  Similarly, hydrodynamic tunneling was observed experimentally quite early on by Eddi \textit{et al.} \cite{EFMC09}, however a rigorous analytical study was conducted much later by Nachbin \textit{et al.} \cite{NMB17}.

On the other hand, phenomena at larger scales that necessitate quantum properties, such as quantum saturable absorbers, transistors, and others, have yet to be studied using the walking droplet system.  In order to make these more complex connections with quantum mechanics it is necessary to develop accurate models of the walking droplet systems.  While there are many complex properties, the main features of the phenomenon are the evolution of the pilot wave, the motion of the droplet, and the interaction between the two via discrete impacts.  On this front, a thorough fluid mechanics study was conducted by Mola\v{c}ek and Bush \cite{MolBush13a, MolBush13b}.  However, it is often necessary to develop more mathematically tractable models that may be simulated using fast numerical techniques.  One such model simplifies the interaction between the droplet and bath by ignoring the vertical dynamics and assuming that on average the droplet is propelled horizontally proportional to the gradient of the Bessel-like wavefield from Oza \textit{et al.} \cite{ORB13}.  Other simplified models are in the form of discrete dynamical systems by Gilet \cite{Gilet14, Gilet16} and Rahman \cite{Rahman18}, which is the focus of this article.

It has been shown that a significant benefit of reduced dynamical systems models is their amenability to rigorous mathematical analysis.  Interesting dynamical behavior such as period doubling and chaos have been observed in detailed hydrodynamic experiments of droplets transitioning from the bouncing to the walking regime \cite{MolBush13a}, however past models were not mathematically tractable enough to rigorously prove these observations.  Through the model of Gilet \cite{Gilet14} it was conjectured that a droplet confined to a rectangular cavity will experience Neimark--Sacker \cite{Neimark, Sacker} bifurcations, which was then proved by Rahman and Blackmore \cite{RahmanBlackmore16}.  Since Gilet's model incorporated both the wavefield and droplet dynamics, the bifurcations, and hence the route to chaos is fairly complex.  Rahman and Blackmore provide quite detailed analyses of new types of homoclinic bifurcations \cite{RJB17} leading to chaos \cite{RahmanBlackmore20} in a class of dynamical systems that include Gilet's model \cite{Gilet14}.  If we consider a simpler domain, such as droplets walking on an annulus, for which we average out the interactions between the droplet and the wavefield, we can derive a 1-D map for the velocity of the droplet.  This makes the analysis far more accessible as shown by Rahman \cite{Rahman18}.  A detailed overview of dynamical models can be found in the review by Rahman and Blackmore \cite{RahmanBlackmoreReview20}.

More broadly, damped-driven systems, such as the hydrodynamic quantum analog, are universal.  Specifically, the energy balance physics is canonical in that it is prevalent in a broader range of  driven-dissipative physical systems that are spatio-temporal in nature, including for example, mode-locked lasers~\cite{haus2000mode,kutz2006mode,ding2009operating,LiWaiKutz10,bale2009transition}, rotating detonation engines (RDEs)~\cite{koch2020mode,koch2020modeling,koch2021multiscale}, and optical communications systems~\cite{agrawal2012fiber,kutz1998hamiltonian,kutz1999dynamics}. In each of these systems, the gain-loss dynamics produces a universal, underlying bifurcation structure.  Specifically, there is an observed period doubling route to chaos as the driving is increased.  This instability cascade has been characterized experimentally in two systems: mode-locked lasers~\cite{bale2009transition} and RDEs~\cite{koch2020mode,koch2020modeling,koch2021multiscale}. The hydrodynamic quantum analog has been experimentally observed~\cite{Tambasco2016,Durey2020c} to produce the same bifurcating cascade  as these systems, showing that that the particle-wave interaction dynamics is a manifestation of canonical damped-driven dynamics in spatially extended systems.

\begin{figure}[htbp]
    \centering
    \stackinset{l}{}{t}{1pt}{\textbf{(a)}}{\includegraphics[height = 0.24\textheight]{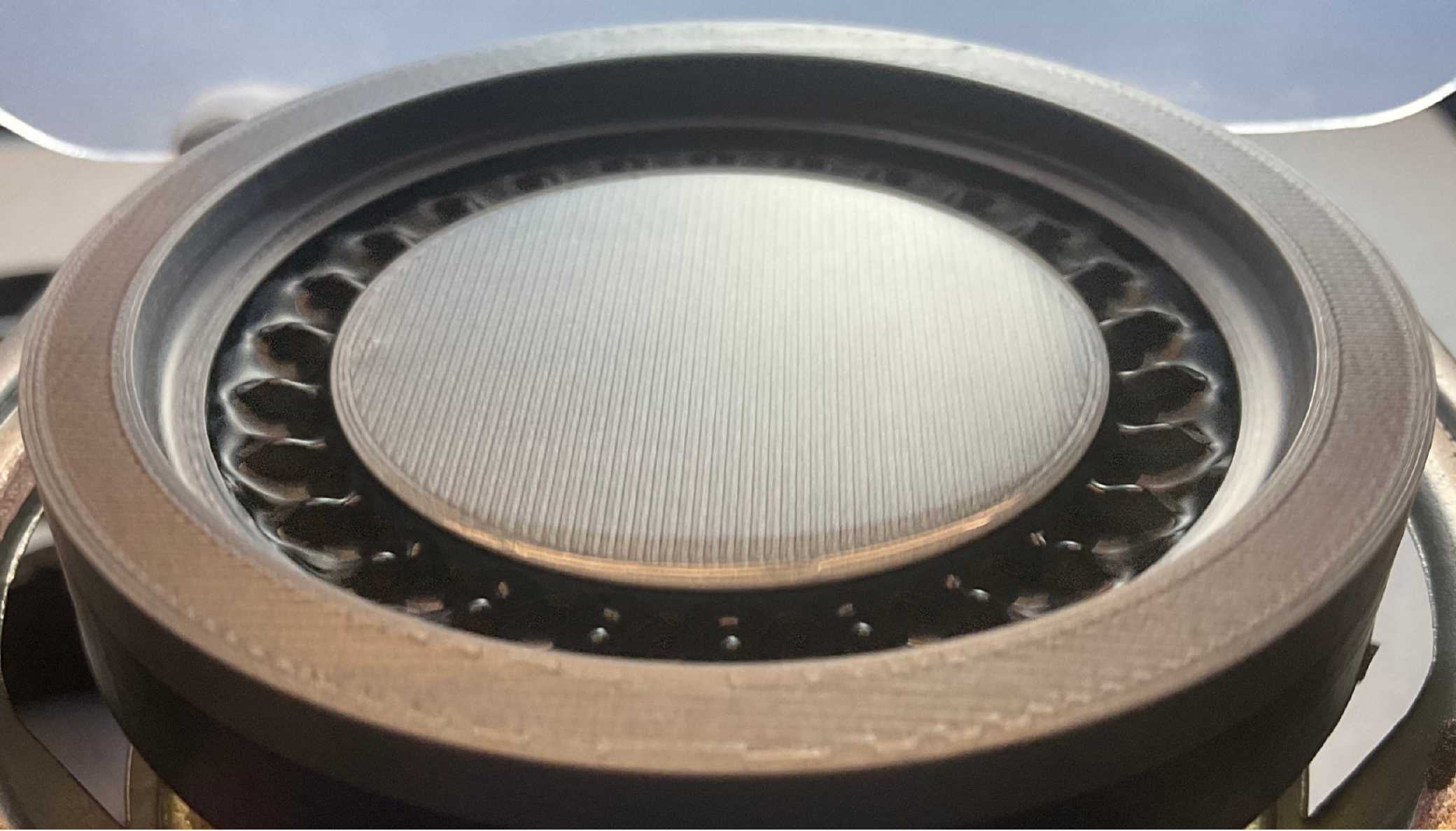}}\,
    \stackinset{l}{}{t}{1pt}{\textbf{\color{white} (b)}}{\includegraphics[height = 0.24\textheight]{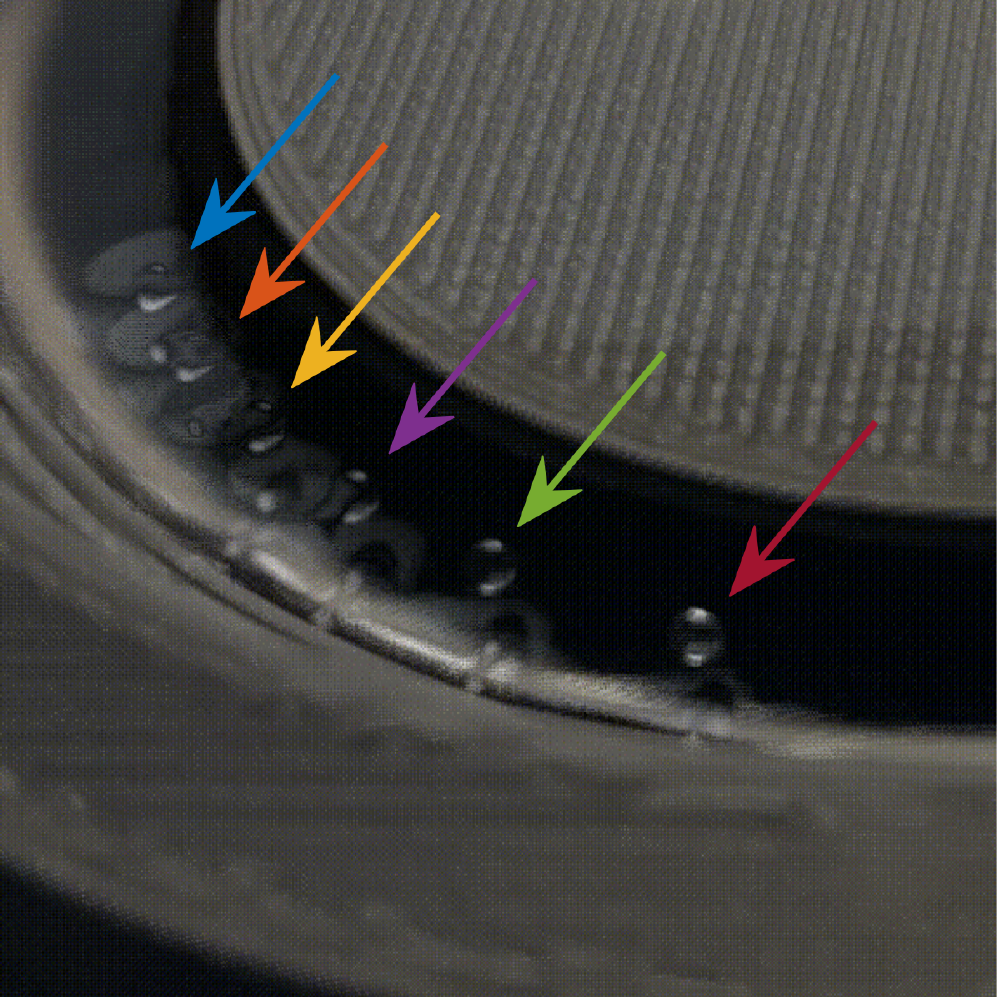}}
    \caption{Experimental backdrop of the theoretical framework.  The bath is vertically forced at $55$ Hz with a maximum acceleration of approximately four times that of gravity.  (a) Faraday waves in an annular cavity.  (b) Sample trajectory of unsteady, but not yet chaotic, walking on the annulus.  The droplet is in the excited walking state with the energy gain (due to the bath) dominating the energy loss (due to the air and deformation).  The chronological progression of the droplet is top left to bottom right:  blue, red, yellow, purple, green, and finally magenta arrows. [Photo credit:  Leviticus James Rhoden]}
    \label{Fig: Experiments}
\end{figure}

The objective of this work is to construct an approximate theoretical framework to model the energy balance in the hydrodynamic quantum analog system.  Much like the theoretical reduction of Li \textit{et al.}~\cite{LiWaiKutz10}, we construct a discrete mapping which approximates the effects of the gain and loss dynamics for walking droplets on an annulus as shown in Fig. \ref{Fig: Experiments}.  The gain and loss mechanisms of droplets are much more complex than the mode-locked laser system considered by Li \textit{et al.}~\cite{LiWaiKutz10} as the droplet-pilot wave interaction is complex, depending both upon the angle of impact of the droplet as well as the phase of the vibrating fluid bath.  Regardless, we are able to construct a mathematically tractable model that highlights the energy balance in the system.  Moreover, these simple energy considerations are sufficient to describe the observed bifurcation sequence of period-doubling to chaos~\cite{Tambasco2016,Durey2020c}.  More broadly, such model reductions reflect the observation of Robert May, who was highly influential in popularizing the logistic map, that simple systems can produce quite complicated behavior~\cite{LogisticMap}.  Indeed, the logistic map is perhaps the simplest system to characterize damp-driven behavior and thus provides a canonical representation of many more complicated systems.

The remainder of the paper is organized as follows:  We begin our discussion by introducing damped-driven systems and their iterative behavior in Sec. \ref{Sec: Damped-driven}.  To illustrate how standard dynamical systems techniques can be adapted to analyze energy gain-loss maps, we present the logistic map in two steps as a damped-driven system, and show that the usual qualitative analysis using cobweb plots also apply to the two step map.  In Sec. \ref{Sec: Kick Model}, we derive an energy gain-loss model of a droplet on an annulus behaving as a kicked rotator \cite{Chirikov1979, KickedRotator, Rahman18}.  Section \ref{Sec: Comparisons} briefly compares the walking droplet system with other damped-driven systems.  Finally, the investigation is concluded in Sec. \ref{Sec: Conclusion} with a discussion on the implications of the analysis and future work.

\section{Damped-Driven Systems as Maps}\label{Sec: Damped-driven}

The theoretical constructs in this work revolve around approximating continuous spatio-temporal systems, typically modeled by partial differential equations, as discrete maps in time whose spatial field have been parametrized in some way.  This construction gives a significant simplification of the model and allows for recurrent iterations to characterize the dynamics.  While rigorous reduction techniques (e.g. \cite{GoodmanMap2008, GoodmanRahman2015}) can be quite complex, we illustrate geometrically that knowing the damping and driving mechanisms is sufficient to derive such a map.

\subsection{Iterating gain and loss}

Our analysis is motivated by the energy gain-loss formulation of mode-locked lasers. Li \textit{et al.}~\cite{LiWaiKutz10} proposed an iterative model that quantifies the interaction of saturable gain and nonlinear loss in a mode-locked laser cavity. The resulting geometric description of the laser dynamics completely characterizes the generic multi-pulsing instability observed in experiments. In this formulation, the mode-locked pulse energy is computed over one round trip. The exact form of the localized mode-locked pulse solution is unimportant, and it is assumed that the effects of chromatic dispersion, self-phase modulation, nonlinear gain, and the bandwidth gain limitations effectively balance each other to form the mode-locked pulse solution~\cite{haus2000mode,kutz2006mode}. Indeed, the fundamental premise of mode-locking is that a localized pulse solution exists for a given balance of dispersive and dissipative effects.  In this manuscript, we aim to use these ideas to derive the energy gain-loss formulation for walking droplets on an annulus.

In the study of Li \textit{et al.}~\cite{LiWaiKutz10} an iteration scheme is constructed by approximating the gain and loss dynamics for each round trip of the cavity as discrete events. The gain dynamics in mode-locked lasers can be modeled in a number of ways, including by a low-order gain saturation differential equation~\cite{LiWaiKutz10}.  The specific form of gain can vary greatly depending on the system being studied, and thus is not the important factor, rather the fact that it produces an input-output relation for the gain 
\begin{equation}
    E_{n+1}^\text{gain} = G(\cdot)
    \label{eq:trans2}
\end{equation}
where $E_{n}$ denotes the energy in the $n\textsuperscript{th}$ round trip of the cavity and where $G(E_n)$ specifies the gain curve was of importance to the study of Li \textit{et al.}~\cite{LiWaiKutz10}, and is also of importance to the present investigation.  

The nonlinear loss in the cavity, i.e., the saturable absorption or saturation fluency curve, is modeled by a transmission function:
\begin{equation}
    E_{n+1}^\text{loss} = T(\cdot)
    \label{eq:trans1}
\end{equation}
The transmission function $T(\cdot)$ can vary significantly from experiment to experiment, especially for very high input energies.  Importantly, this also produces an input-output relationship like (\ref{eq:trans2}).

This separation of the gain and loss gives a composition mapping over a single round trip of the laser cavity
\begin{equation}
    E_{n+1}= T ( G (E_{n}) ) .
\end{equation}
It should be noted that Li \textit{et al.} \cite{LiWaiKutz10} chose the composition $(T\circ G)(\cdot)$, but it would be equivalent to use $(G \circ T)(\cdot)$; i.e., $E_{n+1}= G ( T (E_{n}) )$, where $E_n$ would still represent the energy over a single round trip of the laser cavity.  Fixed points of this map represent energy balanced solutions.  However, because of the mapping itself, a diversity of dynamics can be observed, including the canonical period-doubling bifurcation to chaos~\cite{LiWaiKutz10}.  In the original application of mode-locked lasers, there are numerous theoretical models for accounting for the physical mechanisms for producing nonlinear gains and losses~\cite{haus2000mode,kutz2006mode}, including saturable absorption~\cite{kutz1997mode}, polarization rotation~\cite{spaulding2002nonlinear}, mode-coupling~\cite{proctor2005passive,proctor2005nonlinear}, long-period fiber gratings~\cite{intrachat2003theory}, and phase-sensitive amplification~\cite{kutz2008passive}.  Each of these methods produces their own unique $T(\cdot)$ and $G(\cdot)$, but the dynamics in each case is the same:  a period doubling route to chaos is observed.  Thus the dynamics are observed to be insensitive and robust to the specific form of the transmission functions. 

The goal in what follows is to generalize this framework in order to accommodate it to the hydrodynamic quantum analog system.  Specifically, we develop models from the governing hydrodynamics equations to approximate both $T(\cdot)$ and $G(\cdot)$.  More importantly, we show how the logistic map exhibits many of the universal properties needed for our gain-loss analysis.

\subsection{The logistic map in two steps}

To give a simple example of how a two-step map can be formulated for the gain-loss dynamics and the $T(\cdot)$ and $G(\cdot)$, we consider a reformulation of the logistic map~\cite{LogisticMap}.  Consider the logistic map 
\begin{equation}
    E_{n+1} = rE_n - rE_n^2.
\end{equation}
where $E_n$ is again the energy at the $n\textsuperscript{th}$ iteration.  Suppose we know that most of the population gain occurs in a cycle that can be written as $E_{n+1}^\text{gain} = rE_n^2$.  Then the next cycle needs to include the respective loss that brings the map back to the value of the original logistic map, hence we solve for $E_n$ as a function of $E_{n+1}^\text{gain}$, and plug it back into the logistic map, which yields our loss: $E_{n+1}^\text{loss} = \sqrt{rE_{n+1}^\text{gain}} - E_{n+1}^\text{gain}$.  This yields the full two step system,
\begin{equation}
\begin{split}
E_{n+1}^\text{gain} &= r\left(E^\text{loss}_n\right)^2 = G(E^\text{loss}_n),\\
E_{n+1}^\text{loss} &= \sqrt{rE_{n+1}^\text{gain}} - E_{n+1}^\text{gain} = T(E_{n+1}^\text{gain}) ;
\end{split}.
\label{Eq: Logistic map}
\end{equation}

\begin{figure}[htbp]
\centering
\stackinset{l}{}{t}{}{\textbf{(a)}}{\includegraphics[width = 0.45\textwidth]{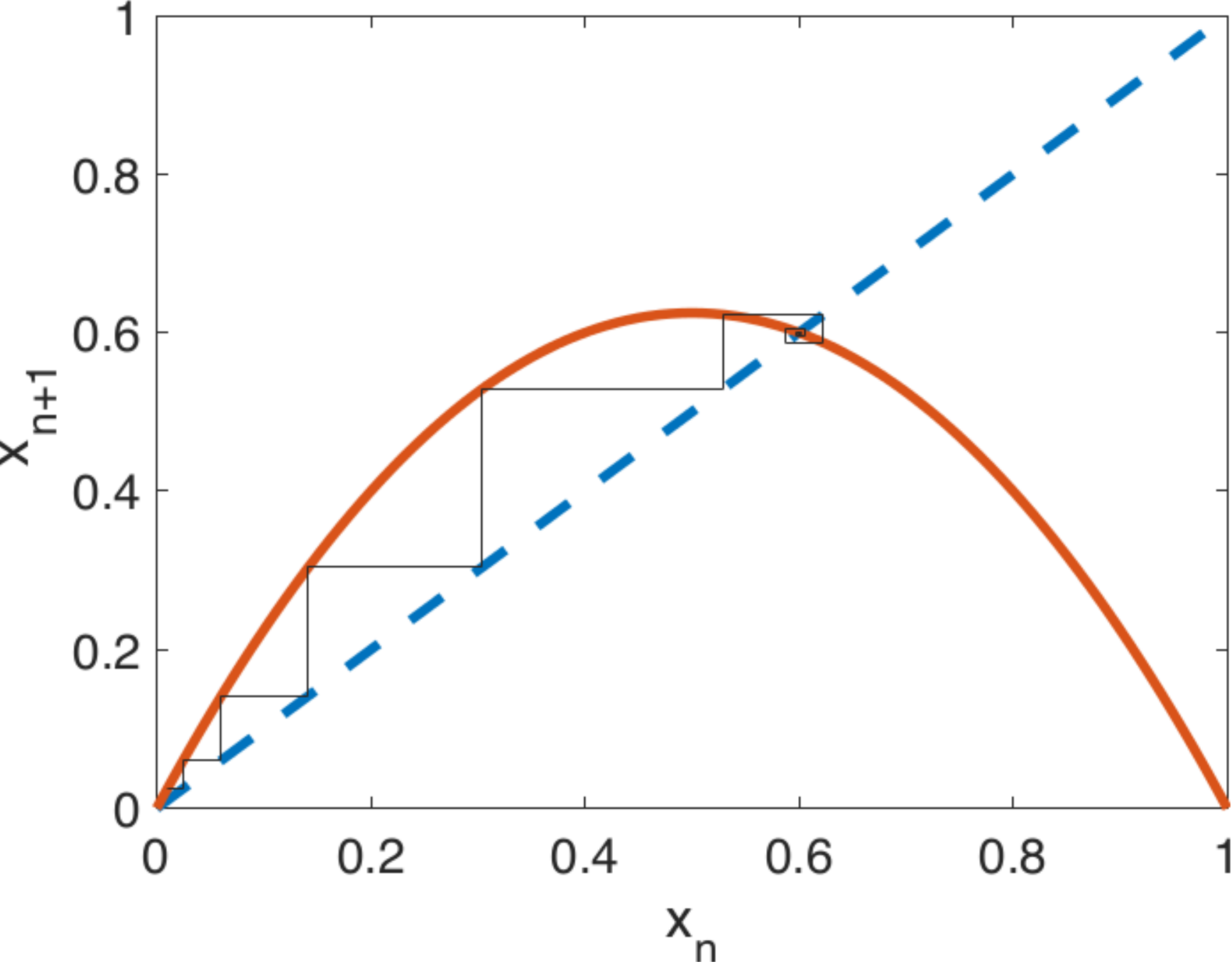}}\qquad
\stackinset{l}{}{t}{}{\textbf{(b)}}{\includegraphics[width = 0.45\textwidth]{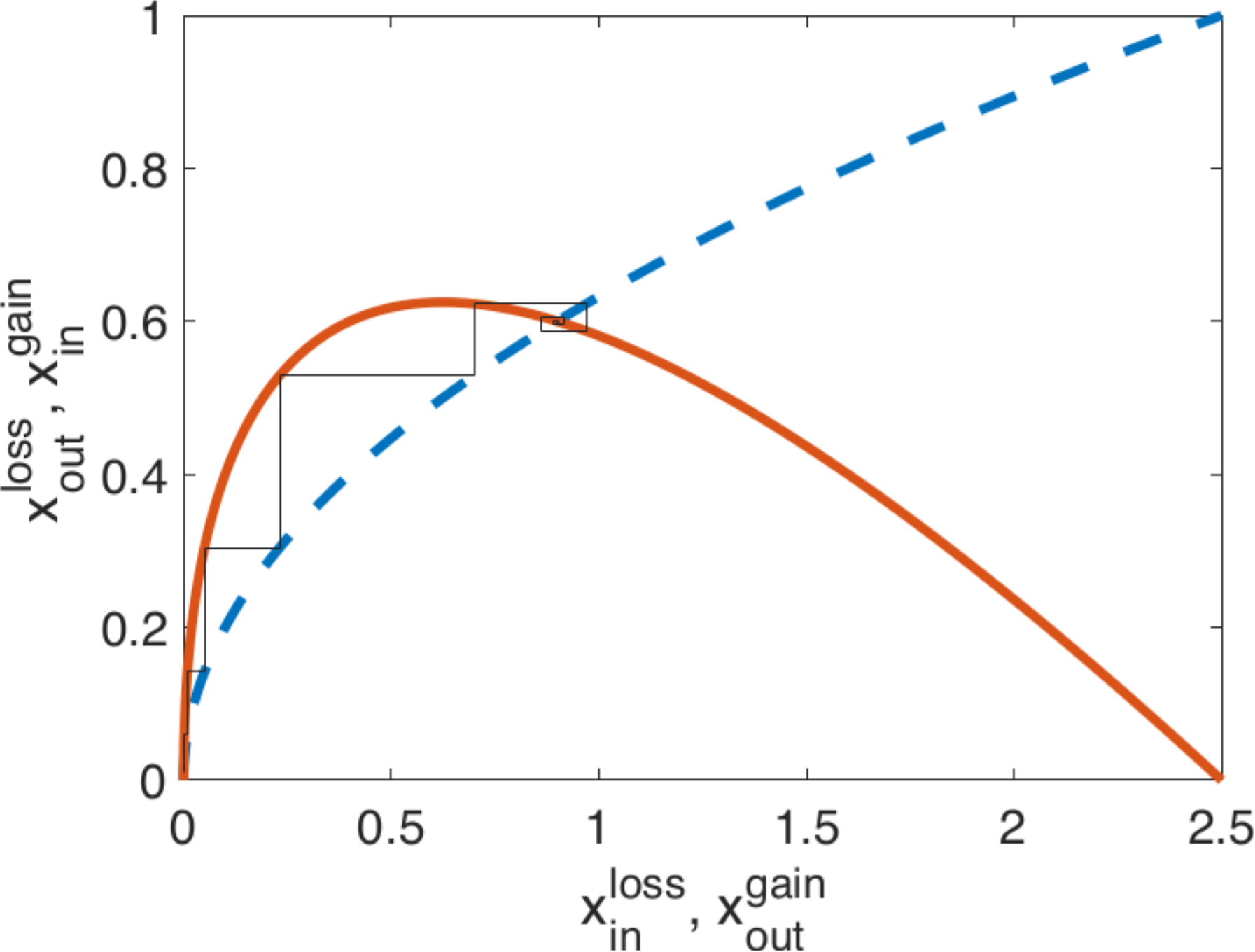}}
\stackinset{l}{}{t}{}{\textbf{(c)}}{\includegraphics[width = 0.45\textwidth]{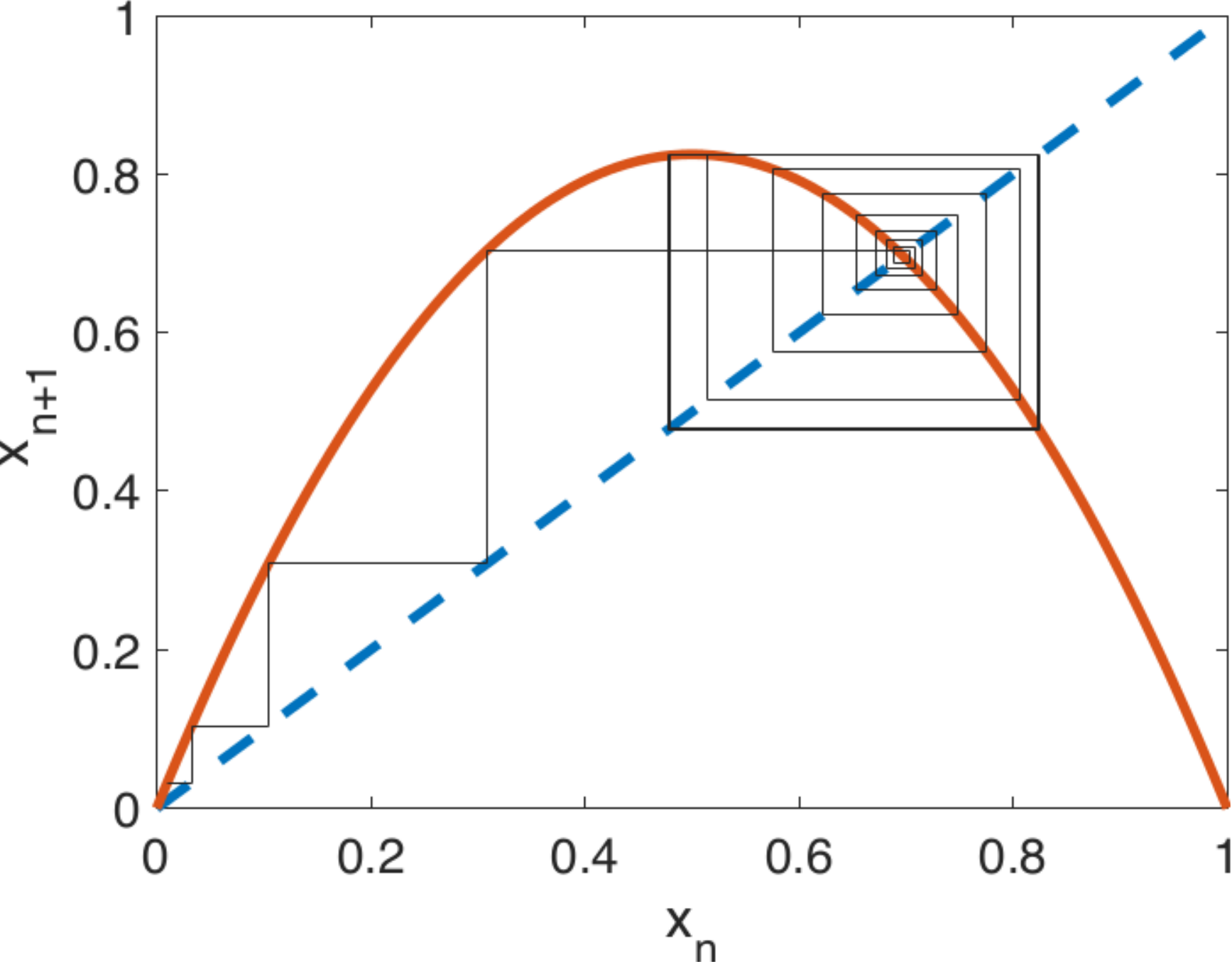}}\qquad
\stackinset{l}{}{t}{}{\textbf{(d)}}{\includegraphics[width = 0.45\textwidth]{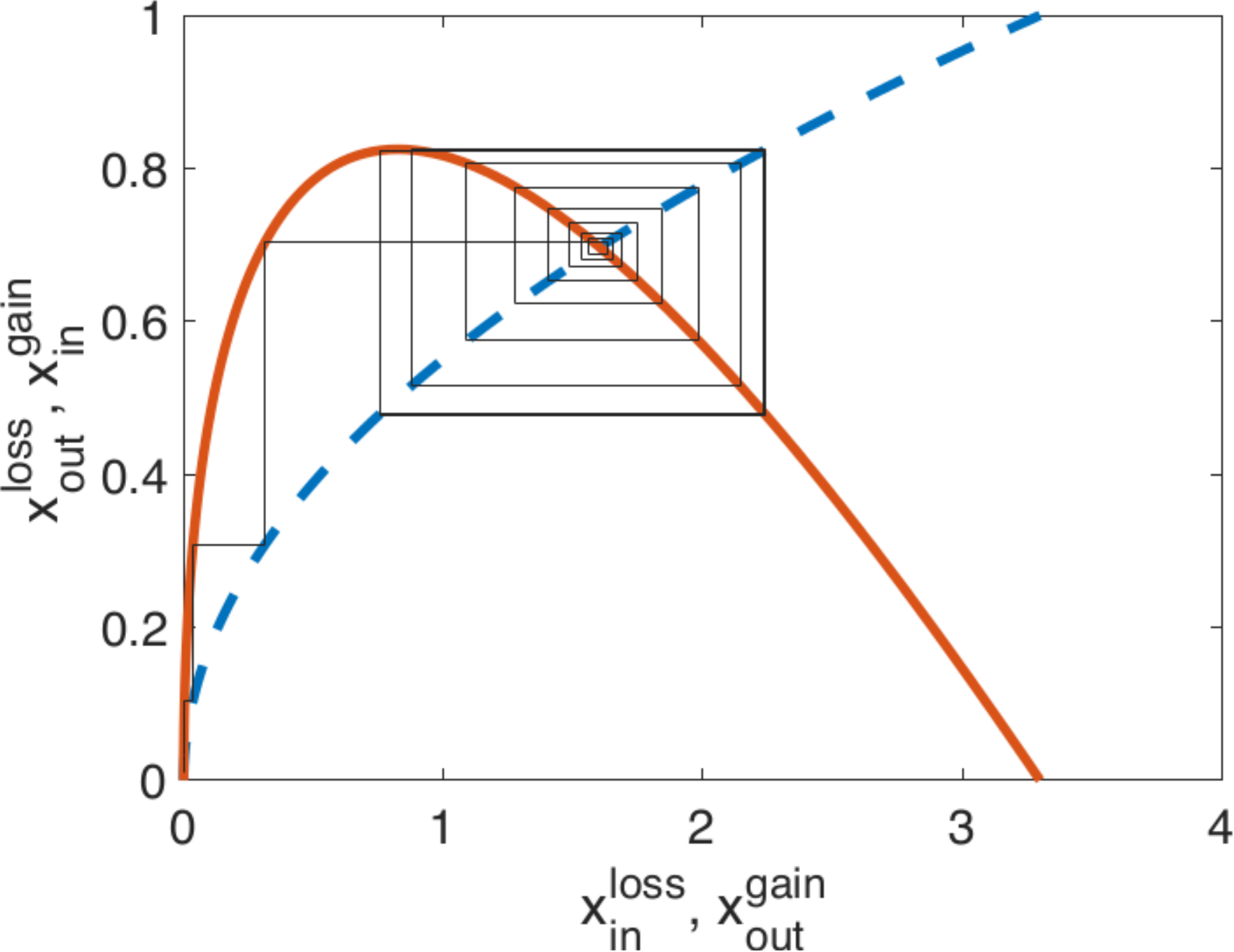}}
\stackinset{l}{}{t}{}{\textbf{(e)}}{\includegraphics[width = 0.45\textwidth]{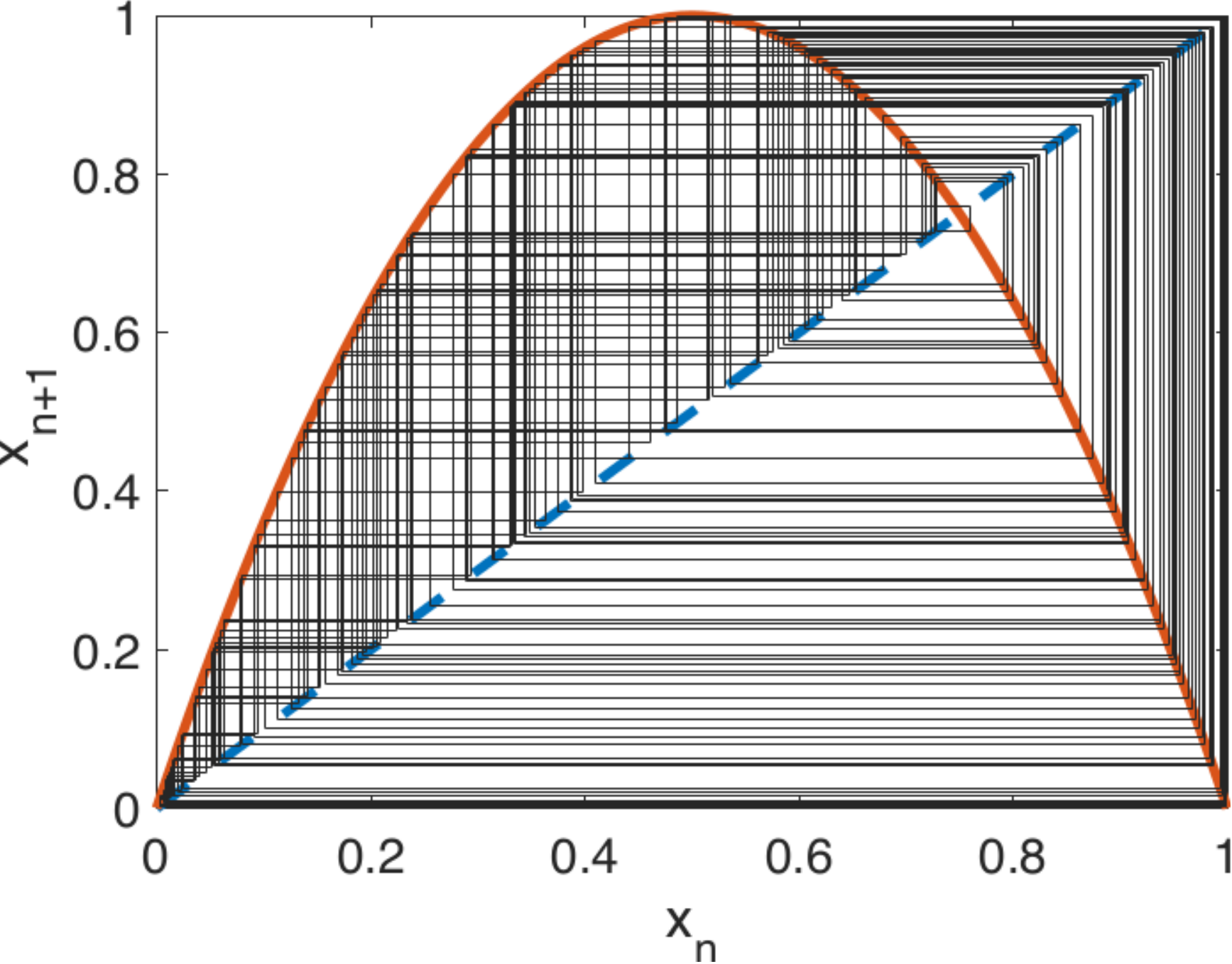}}\qquad
\stackinset{l}{}{t}{}{\textbf{(f)}}{\includegraphics[width = 0.45\textwidth]{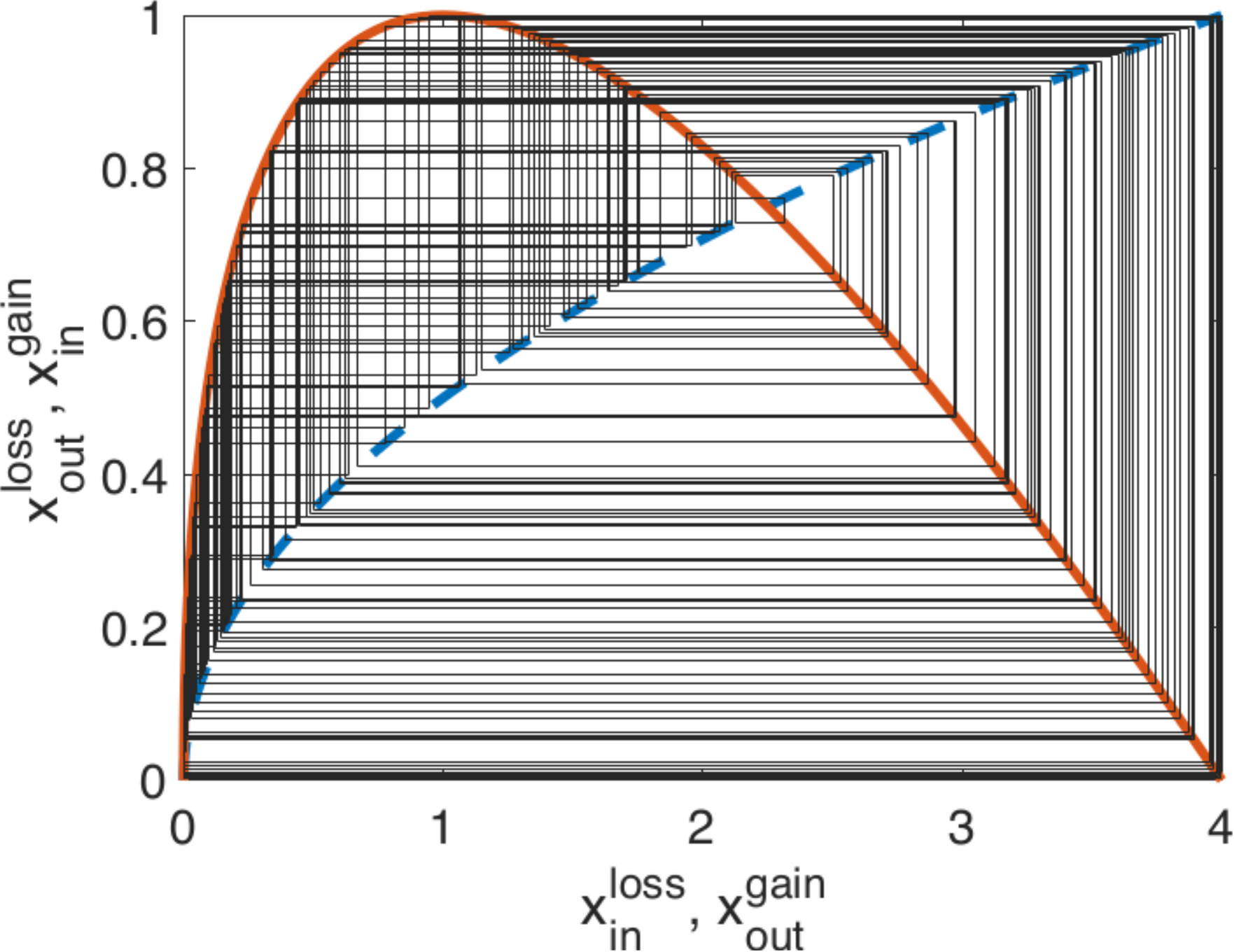}}
\caption{Comparison of the original logistic map \cite{LogisticMap} \textbf{(a, c, e)} and the two step logistic map \eqref{Eq: Logistic map} \textbf{(b, d, f)} for $r = 2.5$ (stable fixed point), $3.3$ (two-cycle), and $4$ (chaotic behavior).  Wide solid (red) curves represent the right hand side of the original Logistic map \textbf{(a, c, e)} and the gain curve of the two step Logistic map \textbf{(b, d, f)}.  Dashed (blue) lines represent the fixed point line for the original Logistic map \textbf{(a, c, e)}, and dashed (blue) curves represent the loss curve of the two step Logistic map \textbf{(b, d, f)}.  The thin dark lines illustrate the cobweb plot.}\label{Fig: Logistic2Step}
\end{figure}

Just as with the original logistic map, we can employ cobweb plots to visualize the dynamics.  In a standard cobweb plot the curve $E_{n+1} = f(E_n)$ and the line $E_{n+1} = E_n$ (i.e., the candidate fixed points) is sketched in Cartesian coordinates.  To create the cobweb, a vertical line from $E_{n+1} = E_n$ to $E_{n+1} = f(E_n)$ is used to iterate to the next value, and a horizontal line from $E_{n+1} = f(E_n)$ to $E_{n+1} = E_n$ is used to iterate the state of the system (for example, $E_{n} \mapsto E_{n+1}$).  To achieve the same results in the two step system \eqref{Eq: Logistic map}, as shown in Li \textit{et al.} \cite{LiWaiKutz10}, we draw a horizontal line from the loss curve to the gain curve (a function with respect to the ordinate) for the gain cycle, and a vertical line from the gain curve to the loss curve (a function with respect to the abscissa) for the loss cycle.  A side-by-side comparison is shown in Fig. \ref{Fig: Logistic2Step}.  The figures provide evidence of the topological equivalence between the iterates of the two maps.

\section{Energy-based formulation of kick-type model}\label{Sec: Kick Model}

In this section we consider the kick model of Rahman \cite{Rahman18}.  Their model describes the averaged changes of velocity for a walker in an annular cavity with bath forcing above the Faraday wave threshold as a discrete dynamical system of the form
\begin{equation}
v_{n+1} = C\left[v_n + K\sin(\omega v_n)e^{-\nu v_n^2}\right].,\\
\label{Eq: KickModel}
\end{equation}
where $v_n$ is the velocity after the $n^\textsuperscript{th}$ impact, $C \ll 1$ is the damping factor, and $K = \frac{-\pi e^{\nu\pi^2}}{\sin(\pi\omega)}$ is the kick strength derived analytically to satisfy functional properties on the annulus \cite{Rahman18}.  Further, $\omega = 31/2$ and $\nu = \omega^2/8.4\pi^2$ are the wavenumber and damping parameter from the experiments of Filoux \textit{et al.} \cite{FHV15}.  In their study, Rahman \cite{Rahman18}, rigorously analyzed the model and used standard Dynamical Systems tools, such as cobweb plots, to illustrate the analysis.  This gives us an opportunity to compare the cobweb plots from the energy formulation with the standard Dynamical Systems formulation, which may not always be possible to derive for all systems, such as the systems of Li \textit{et al.} \cite{LiWaiKutz10} and Koch \cite{koch2020modeling}.

First let us write the recurrence relation in terms of the nondimensional kinetic energy,
\begin{equation*}
v_{n+1}^2 = C^2v_n^2 + C^2K^2\sin^2(\omega v_n)e^{-2\nu v_n^2} + C^2Kv_n\sin(\omega v_n)e^{-\nu v_n^2},
\end{equation*}
which can be written as
\begin{equation}
E_\text{out} = C^2E_\text{in} + C^2K^2\sin^2(\omega \sqrt{E_\text{in}})e^{-2\nu E_\text{in}} + C^2K\sqrt{E_\text{in}}\sin(\omega \sqrt{E_\text{in}})e^{-\nu E_\text{in}}.
\label{Eq: Kick Energy}
\end{equation}
In this setup, all the energy loss will come from the damping factor, and therefore we can write $E_\text{out}^\text{loss} = C^2E_\text{in}$.  Solving for $E_\text{in}$ and plugging into \eqref{Eq: Kick Energy} gives us
\begin{equation*}
E_\text{out} = E_\text{out}^\text{loss} + C^2K^2\sin^2(\omega \sqrt{E_\text{out}^\text{loss}}/C)e^{-2\nu E_\text{out}^\text{loss}/C^2} + CK\sqrt{E_\text{out}^\text{loss}}\sin(\omega \sqrt{E_\text{out}^\text{loss}}/C)e^{-\nu E_\text{out}^\text{loss}/C^2}.
\end{equation*}
We define $E_\text{out}^\text{loss}$ as the energy of the system after incurring losses, and $E_\text{out}^\text{gain}$ as the energy of the system after gaining energy.  Further, since this is a recurrence relation as illustrated in Fig. \ref{Fig: Schematic}, $E_\text{in}^\text{gain} = E_\text{out}^\text{loss}$ and $E_\text{in}^\text{loss}= E_\text{out}^\text{gain}$; that is, after applying the gain mechanism the energy increases from $E_\text{in}^\text{gain}$ to $E_\text{out}^\text{gain}$ and this energy is fed into the loss mechanism that decreases the energy from $E_\text{in}^\text{loss}$ to $E_\text{out}^\text{loss}$, which is then fed back into the gain mechanism, then we can write the loss and gain functions as
\begin{subequations}
\begin{align}
E_\text{out}^\text{loss} &= C^2E_\text{in}^\text{loss},\label{Eq: Kick Loss curve}\\
E_\text{out}^\text{gain} &= E_\text{in}^\text{gain} + C^2K^2\sin^2(\omega \sqrt{E_\text{in}^\text{gain}}/C)e^{-2\nu E_\text{in}^\text{gain}/C^2}\nonumber\\
&\hspace{42pt} + CK\sqrt{E_\text{in}^\text{gain}}\sin(\omega \sqrt{E_\text{in}^\text{gain}}/C)e^{-\nu E_\text{in}^\text{gain}/C^2}\label{Eq: Kick Gain curve}
\end{align}
\label{Eq: Kick Loss-Gain curves}
\end{subequations}
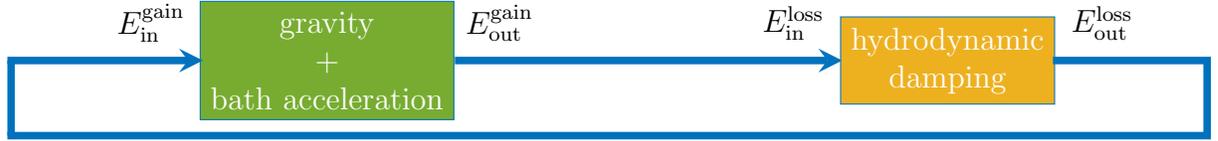
\begin{figure}[htbp]
\centering
\begin{tikzpicture}
\color{blue} 
\node[draw, align = center, fill = green] at (-1.75,0) {\color{white} gravity\\ \color{white}+\\ \color{white}bath acceleration};
\draw[-stealth, line width = 1mm] (-6,0) -- (-3.4,0);
\draw (-3.5 ,0.5) node[left] {\color{black} $E_\text{in}^\text{gain}$};
\draw (-0.05,0.5) node[right] {\color{black} $E_\text{out}^\text{gain}$};
\draw[-stealth, line width = 1mm] (-0.05,0) -- (5.1,0);
\draw (5,0.5) node[left] {\color{black} $E_\text{in}^\text{loss}$};
\node[draw, align = center, fill = yellow] at (6.5,0) {\color{white} hydrodynamic\\ \color{white} damping};
\draw (8,0.5) node[right] {\color{black} $E_\text{out}^\text{loss}$};
\draw[line width = 1mm] (7.9,0) -- (10,0);
\draw[line width = 1mm] (9.95, 0) -- (9.95, -1);
\draw[line width = 1mm] (10, -1) -- (-6,-1);
\draw[line width = 1mm] (-5.95,-1) -- (-5.95,0);
\end{tikzpicture}
\caption{Schematic of the interaction between the driving phase and damping phase as a recurrence relation.  The energy gain in the model comes from the ``kick'' term, $K$ (similar to the standard map), which implicitly captures the velocity of the droplet relative to the bath, and the energy loss comes from the damping factor $C$, which captures the averaged damping due to all hydrodynamic effects.}\label{Fig: Schematic}
\end{figure}

\begin{figure}[htbp]
    \centering
    \stackinset{l}{}{t}{-12pt}{\textbf{(a.i)}}{\includegraphics[width = 0.284\textwidth]{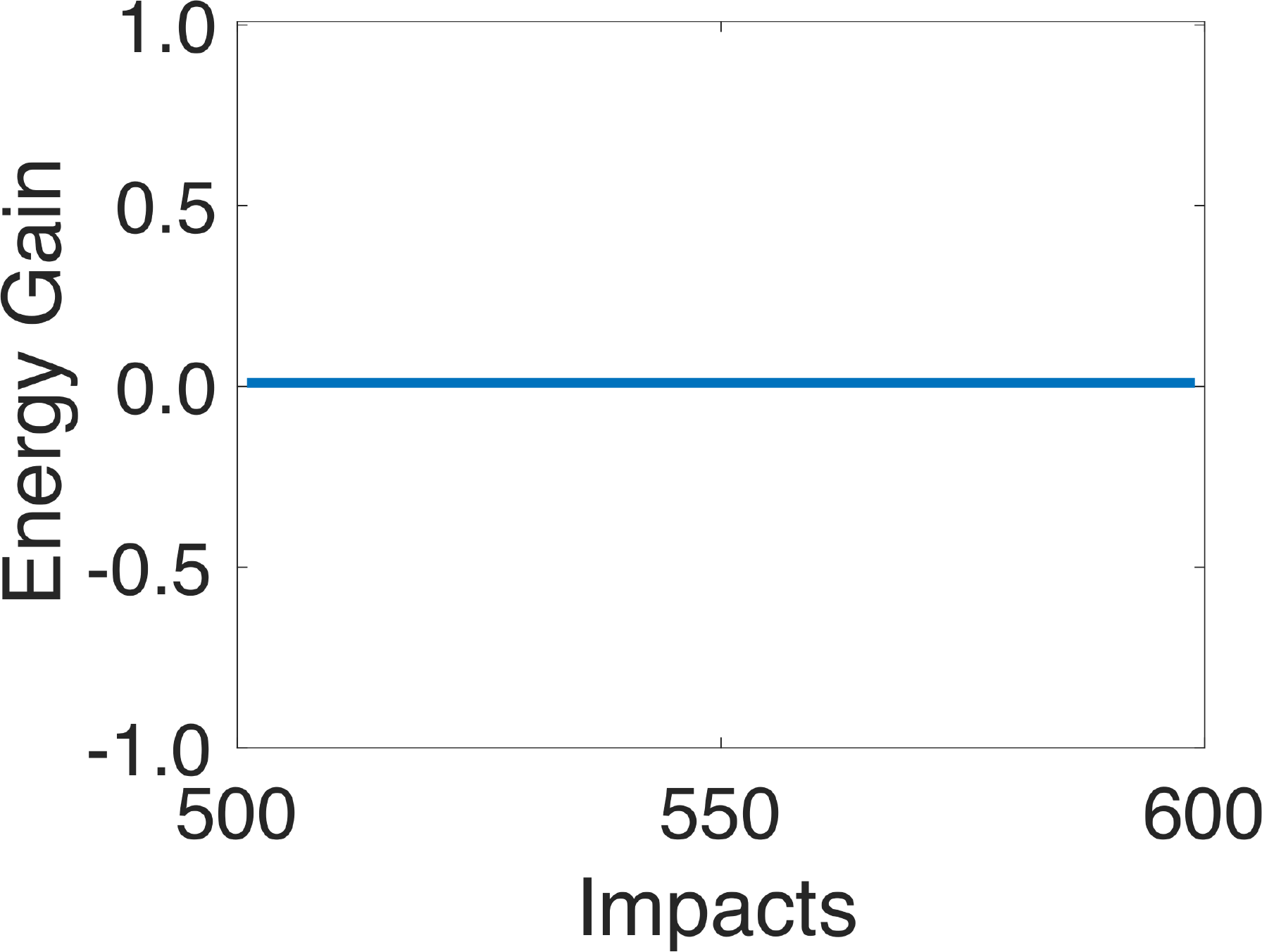}}\qquad
    \stackinset{l}{}{t}{-12pt}{\textbf{(b.i)}}{\includegraphics[width = 0.29\textwidth]{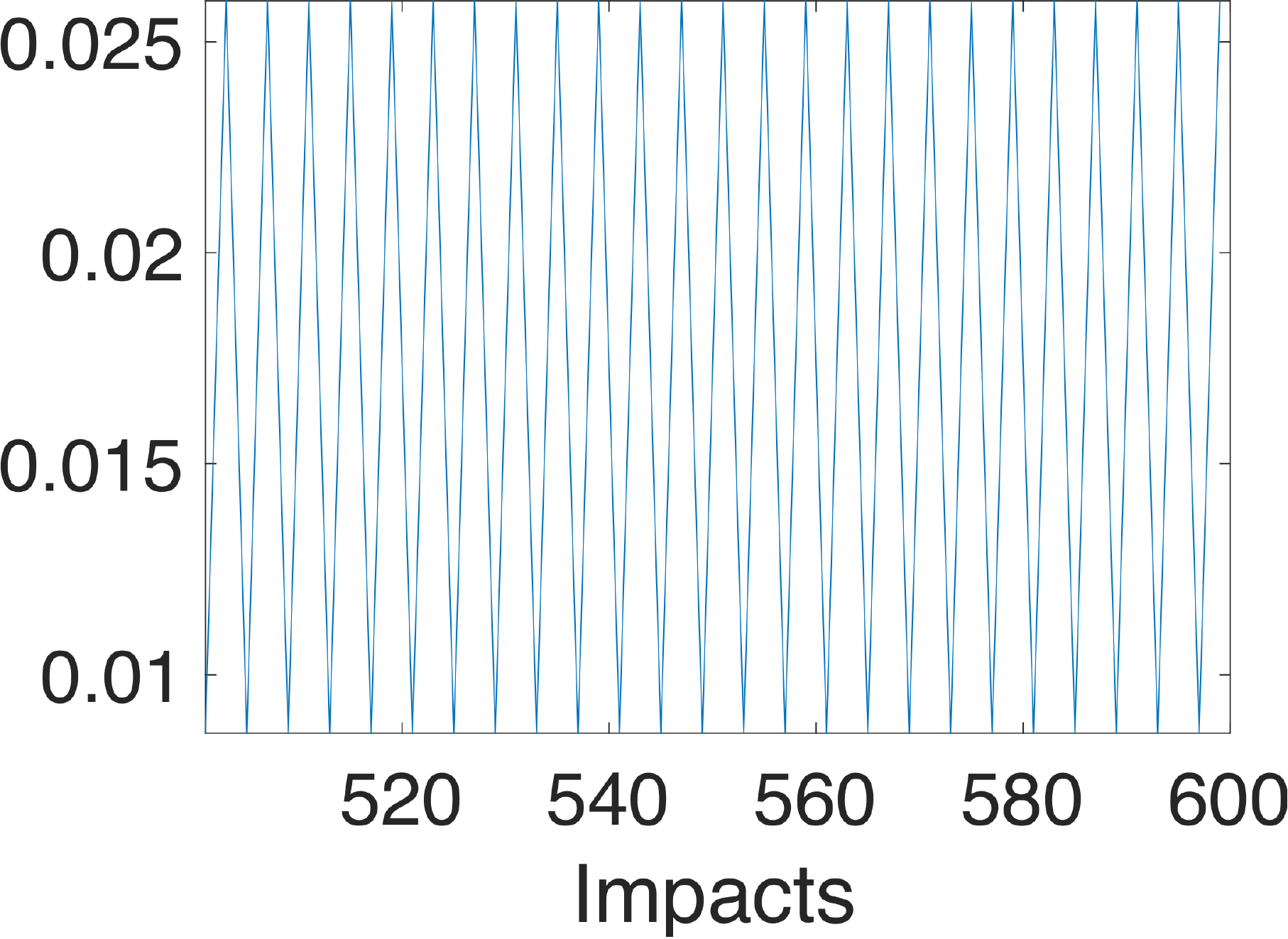}}\qquad
    \stackinset{l}{}{t}{-12pt}{\textbf{(c.i)}}{\includegraphics[width = 0.29\textwidth]{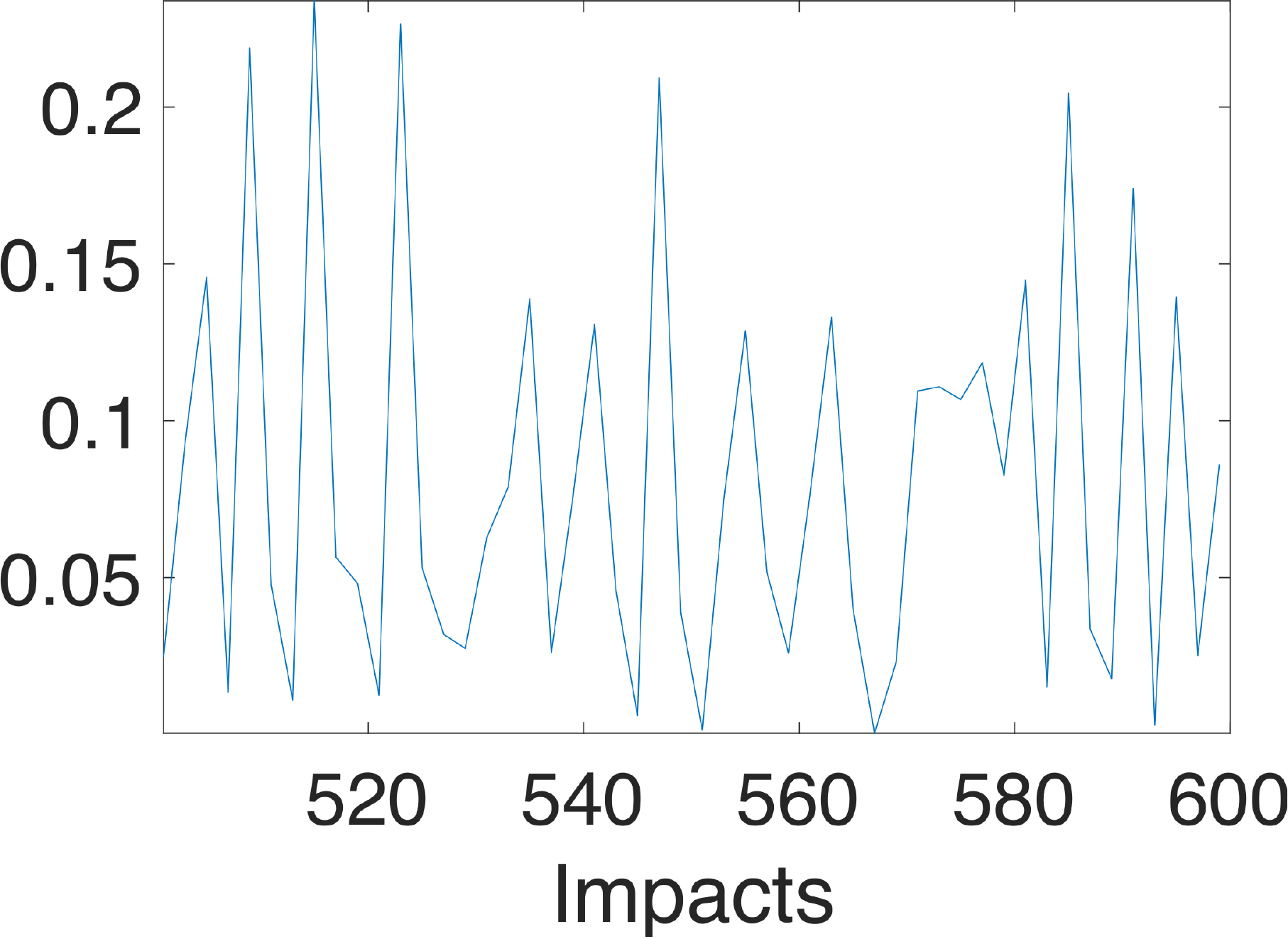}}
    
    \bigskip
    
    \stackinset{l}{}{t}{-12pt}{\textbf{(a.ii)}}{\includegraphics[width = 0.29\textwidth]{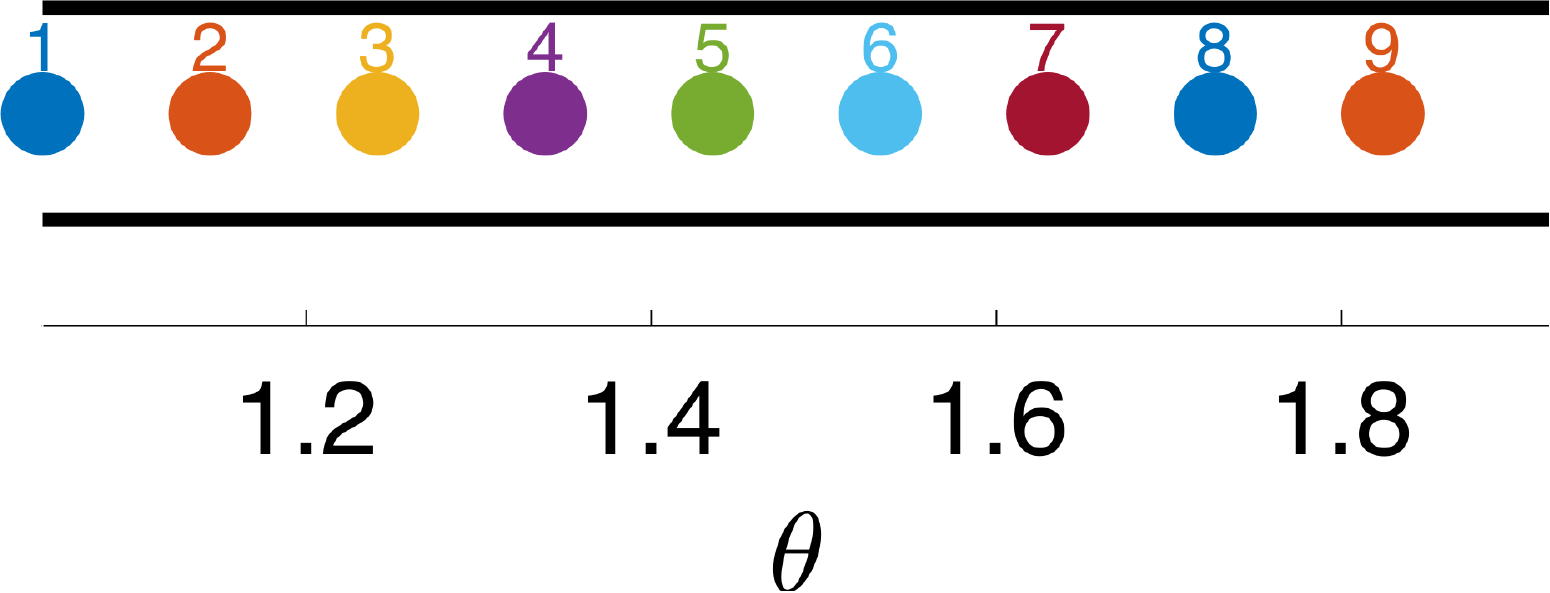}}\qquad
    \stackinset{l}{}{t}{-12pt}{\textbf{(b.ii)}}{\includegraphics[width = 0.29\textwidth]{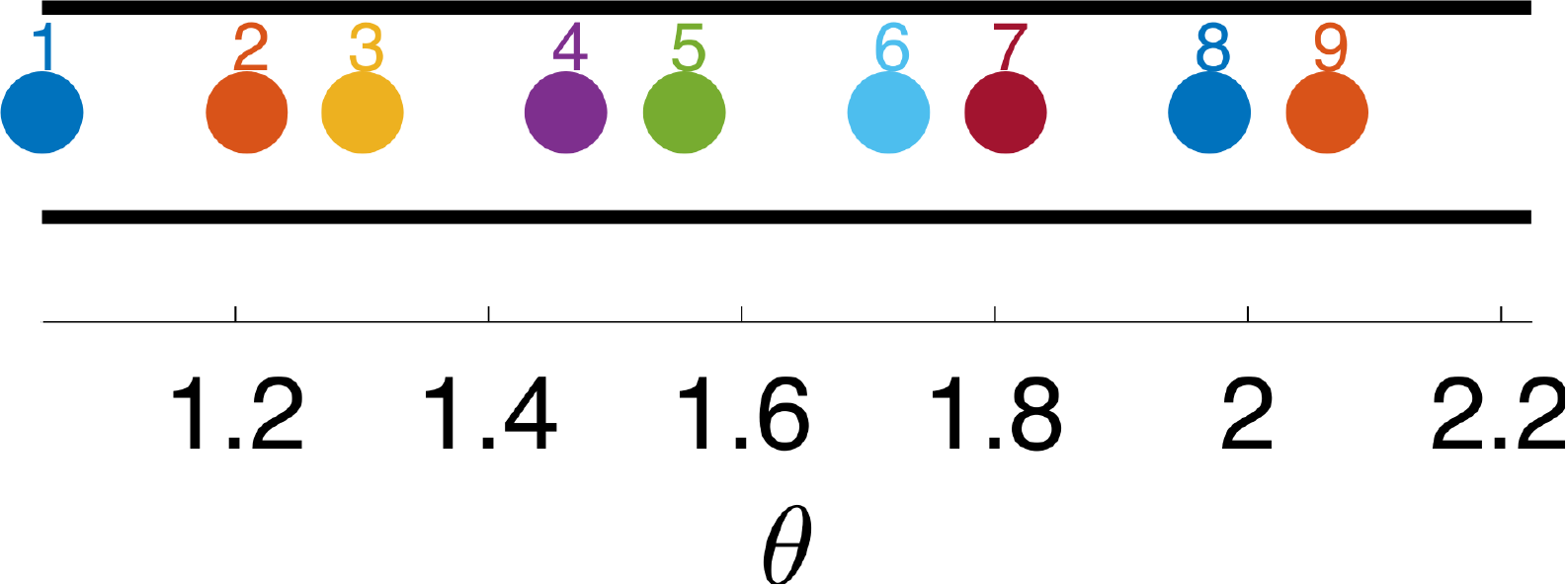}}\qquad
    \stackinset{l}{}{t}{-12pt}{\textbf{(c.ii)}}{\includegraphics[width = 0.29\textwidth]{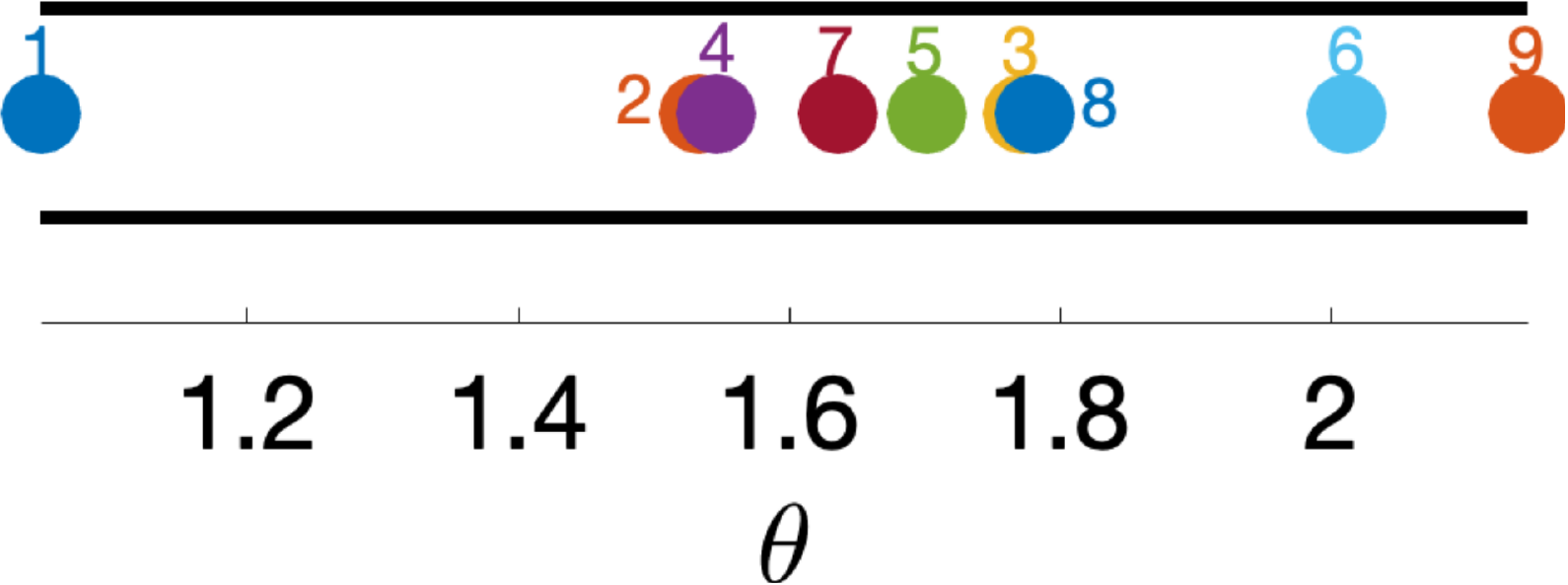}}
    \caption{Transition from the steady walking phase with $C = 1/10K$ (a) to the two-step walking with $C = 1/6K$ (b), and finally to chaotic walking with $C = 1/2K$ (c).  (a.i, b.i, c.i) Time series plots for the energy gain, $E_\text{out}^\text{gain}$.  (a.ii, b.ii, c.ii) Sample trajectories, from simulating the standard map-like model of \cite{Rahman18}, of nine subsequent impacts for the walker on an annulus (straightened out for visual purposes).}
    \label{Fig: Droplet Trajectory}
\end{figure}

The system is sensitive to the damping factor $C$ and the energy gained or lost will depend on the relation between $C$ and $K$, which we illustrate in Fig .  \ref{Fig: Droplet Trajectory}.  When the energy gained and lost are balanced, for example at $C = 1/10K$, the walker walks at a steady velocity as shown in the velocity map \cite{Rahman18}, which has also been observed experimentally \cite{FHV15}.  It should be noted that the energy after one round trip of the gain-loss mechanism could also be balanced if the speed remained constant but the direction switched, however we know from experiments (such as \cite{FHV15}) that the droplet tends to move in the same direction.  The droplet does not switch direction until the velocity destabilizes as shown theoretically in \cite{Rahman18}.

The droplet then transitions to chaotic walking.  While there have been experimental studies on exotic states of walking dynamics below the Faraday wave threshold, such as that by Wind-Willassen \textit{et al.} \cite{WMHB13}, the difficulty of precise measurements at the point of bifurcation makes experiments above the Faraday wave threshold prohibitive in this parameter regime.  In contrast, from the model in the present study and the velocity map of \cite{Rahman18} we observe a period doubling bifurcation, and when we increase $C$ to $C = 1/6K$ we observe a two cycle in the energy domain, which manifests itself as a two-step jump where the droplet bounces a short distance followed by a larger distance.  Finally, the droplet trajectory becomes chaotic, for example at $C = 1/2K$, where the droplet motion is quite disordered and switches between clockwise and counter-clockwise directions on the annulus.

\begin{figure}[htbp]
\centering
\includegraphics[width = 0.9\textwidth]{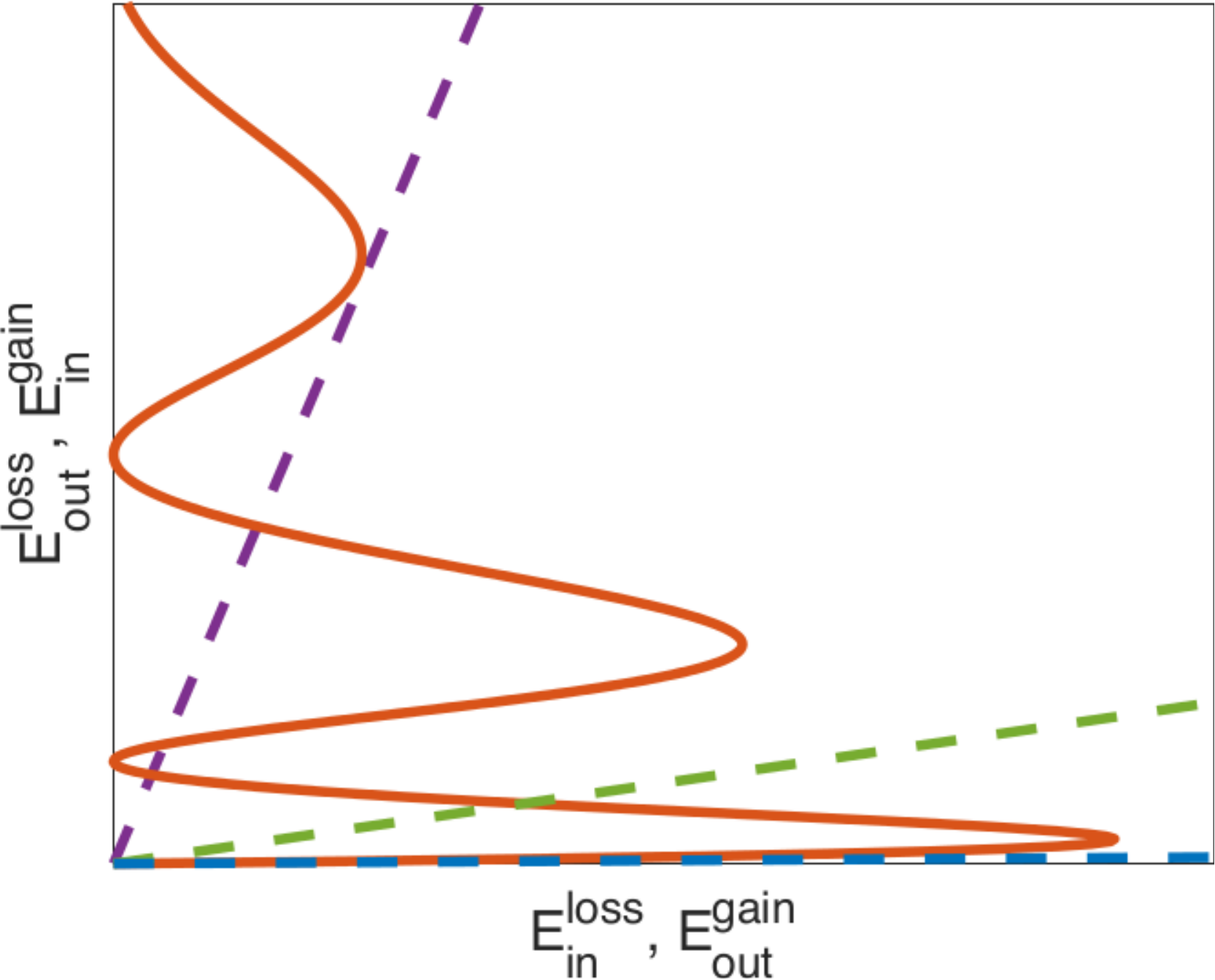}
\caption{Relation between the loss (dashed lines) \eqref{Eq: Kick Loss curve} and gain (solid red) curves \eqref{Eq: Kick Gain curve} at different scales (hence the lack of axis values).  From smallest to largest slope, the loss curves represent \eqref{Eq: Kick Loss curve} with $C = 1/20K$ (blue), $C = 1/4K$ (green), and $C = 1/K$ (purple).  For $C = 1/20K$ there are no non-trivial fixed points, for $C = 1/4K$ there exists a single non-trivial fixed point, and at $C = 1/K$ there are several non-trivial fixed points.}\label{Fig: Kick Gain Loss}
\end{figure}

In this setup, increasing the bifurcation parameter $C$ will give us the same bifurcations that were analytically proved to exist in \cite{Rahman18} as shown in Fig. \ref{Fig: Kick Gain Loss}.  In the figure the system bifurcates from having only a trivial fixed point to having two fixed points (the other being non-trivial) as the parameter is increased from $C = 1/20K$ to $C = 1/4K$.  Further bifurcations occur as we increase the parameter from $C = 1/4K$ to $C = 1/K$ where the system goes through an increasing number of intersections including both transverse and tangential.

\begin{figure}[htbp]
\centering
\stackinset{l}{}{t}{}{\textbf{(a)}}{\includegraphics[width = 0.45\textwidth]{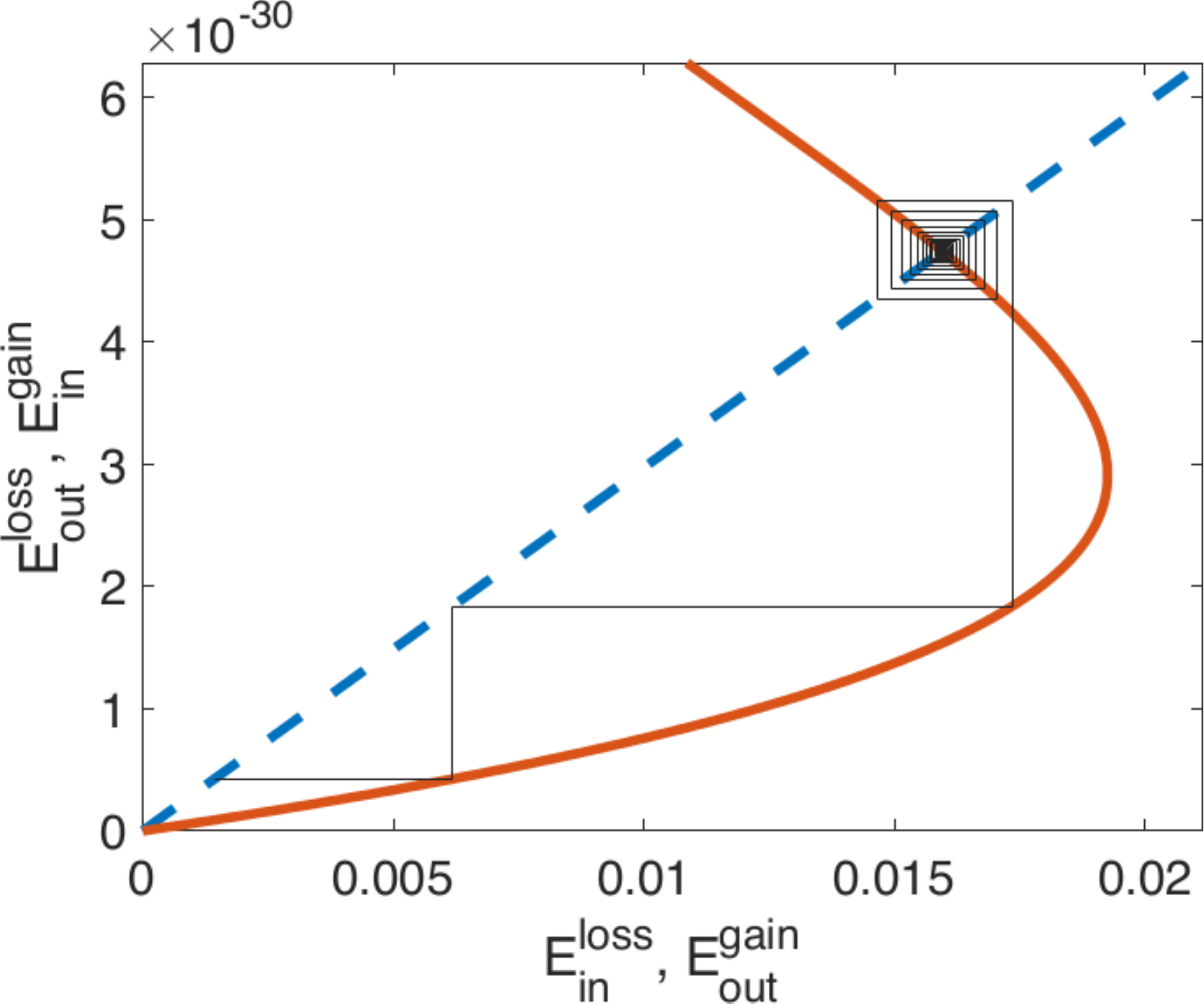}}\qquad
\stackinset{l}{}{t}{}{\textbf{(b)}}{\includegraphics[width = 0.45\textwidth]{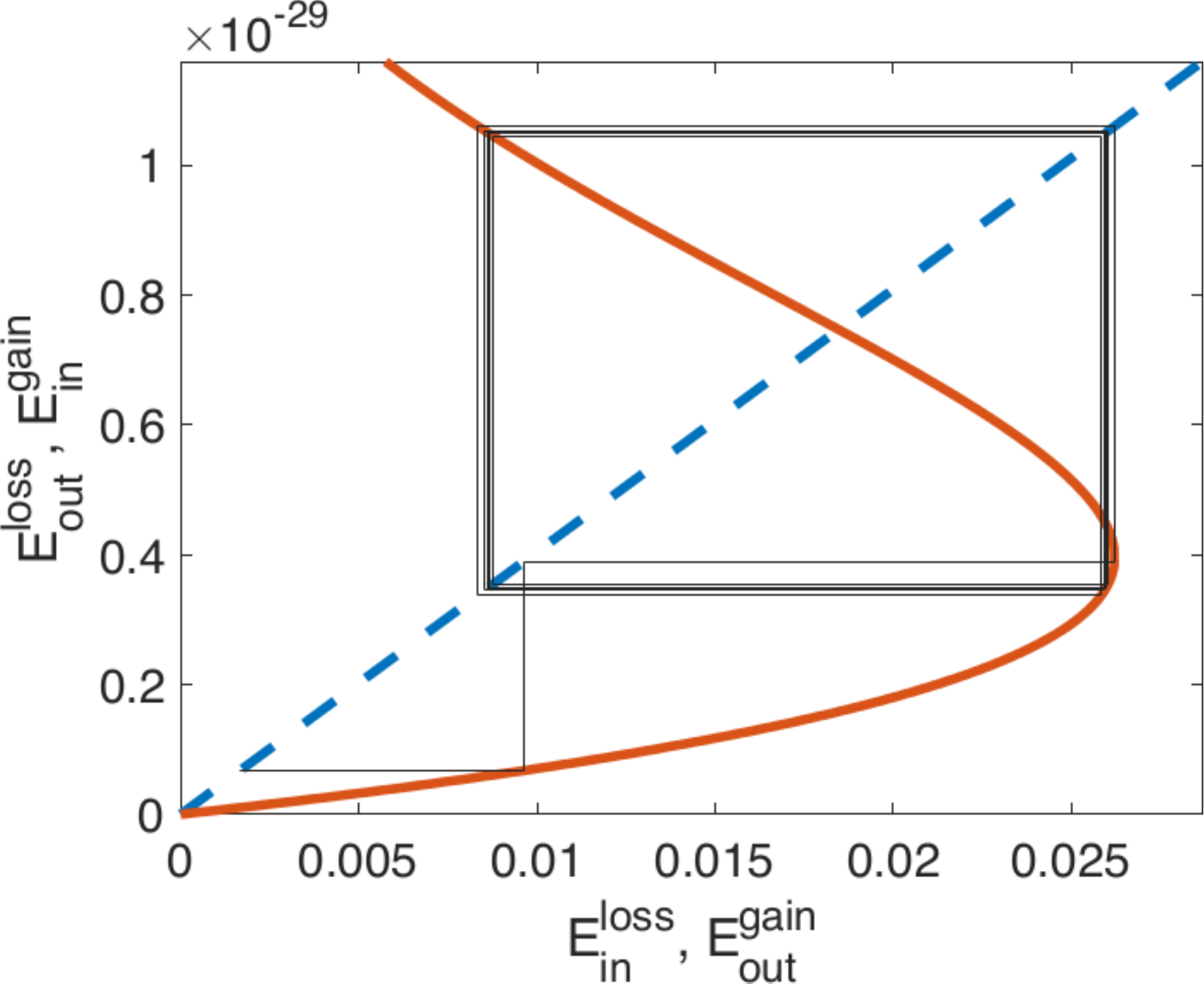}}
\stackinset{l}{}{t}{}{\textbf{(c)}}{\includegraphics[width = 0.45\textwidth]{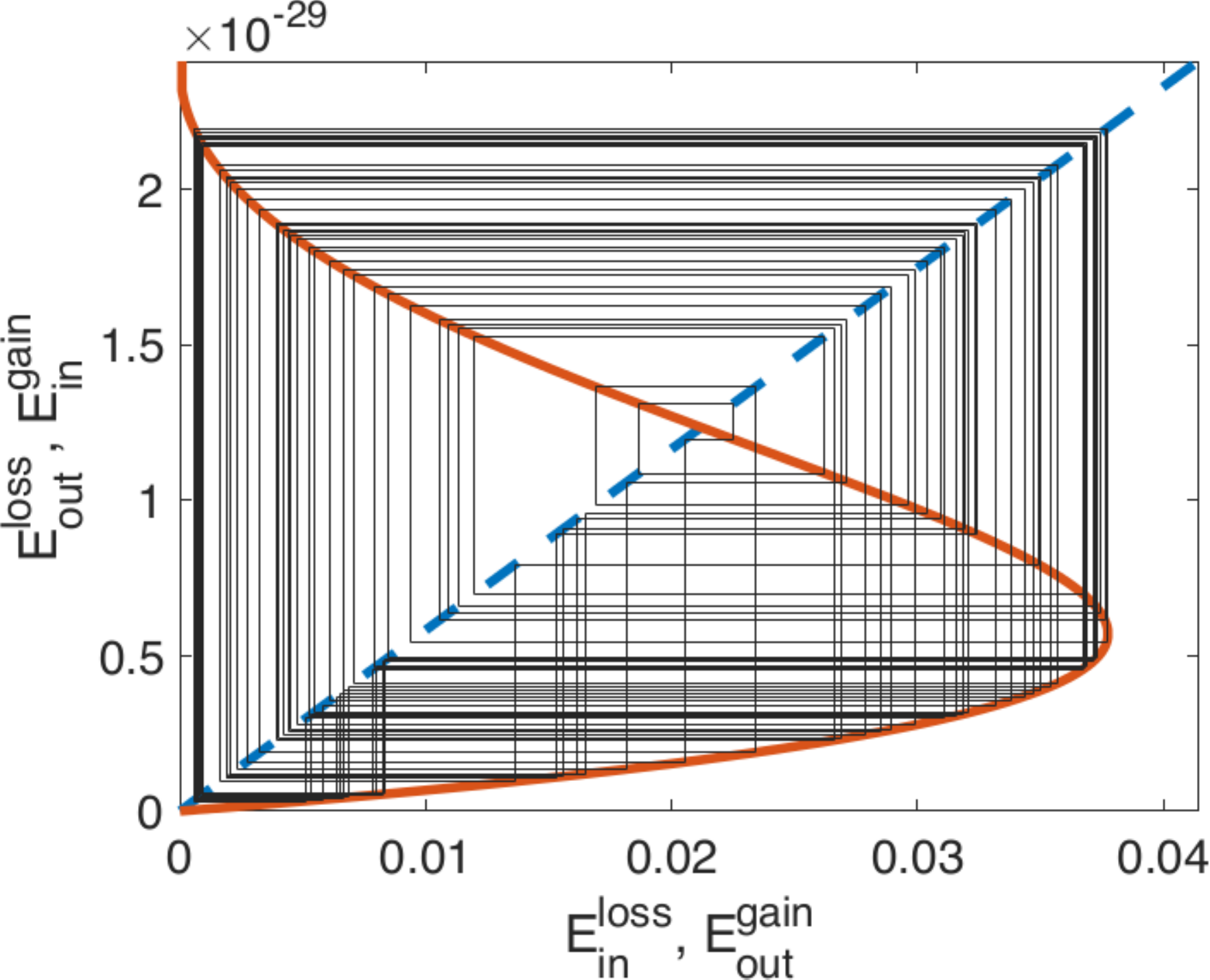}}\qquad
\stackinset{l}{}{t}{}{\textbf{(d)}}{\includegraphics[width = 0.45\textwidth]{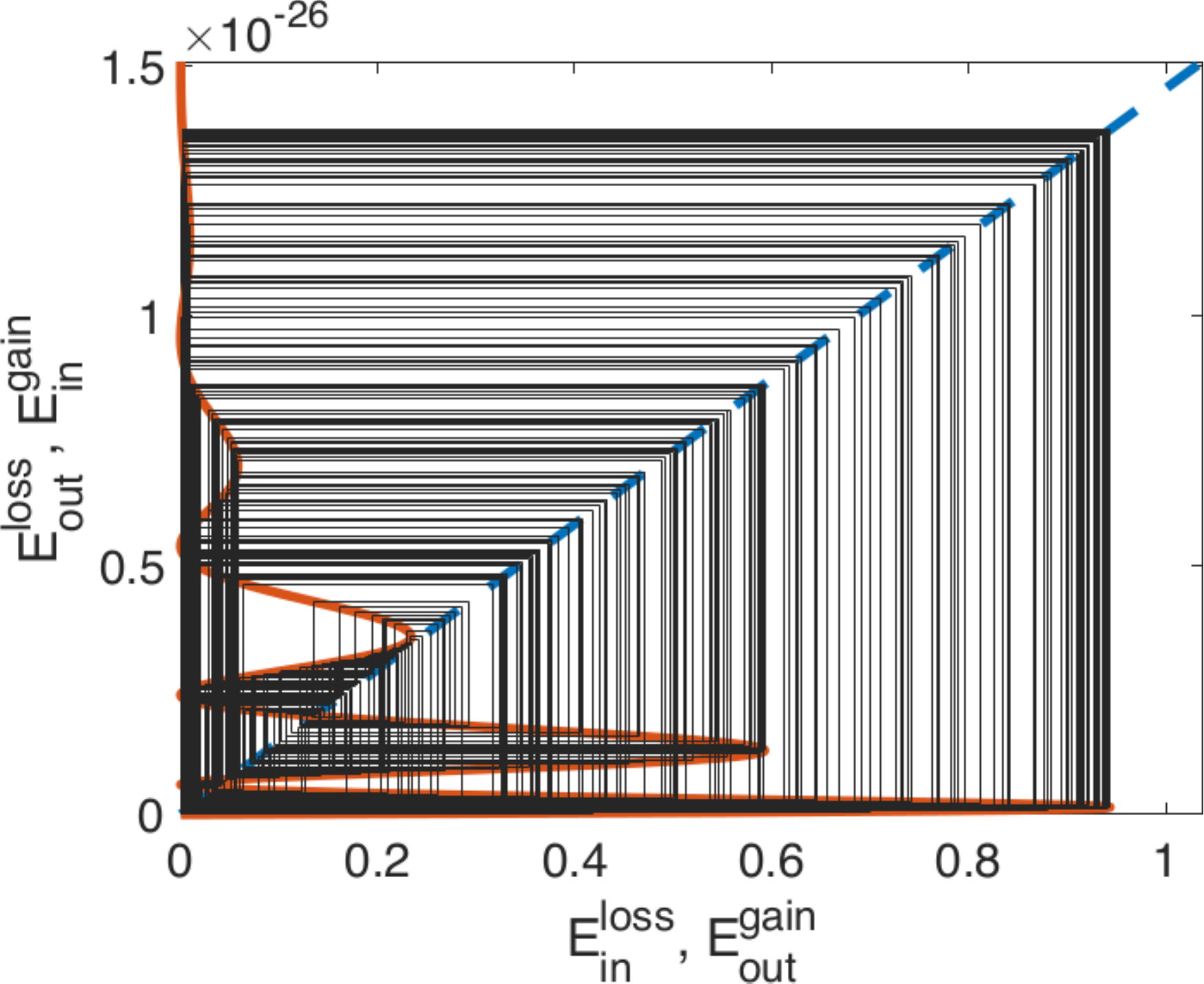}}
\caption{Cobweb plots illustrating the route to chaos.  \textbf{(a)} Convergence to non-trivial fixed point for $C = 1/7K$.  \textbf{(b)} Convergence to a low-period orbit for $C = 1/6K$.  \textbf{(c)} Evidence of seemingly random orbits of arbitrarily large period emerge for $C = 1/5K$.  \textbf{(d)}  Full-scale chaos for $C = 1/K$ as proved in \cite{Rahman18}.}\label{Fig: Kick Transition}
\end{figure}

The dynamics of the damped-driven system of \eqref{Eq: Kick Loss-Gain curves} is further illustrated in the form of cobweb plots.  In \cite{Rahman18}, Rahman illustrated the route to chaos using $C = 1/20K$, $C = 1/4K$, and $C = 1/K$ to complement the proofs, however it is shown that $C = 1/4K$ contains a 3-cycle window, and hence based on the proofs of \cite{Rahman18, LiYorke75} we expect seemingly random trajectories to start appearing for $C < 1/4K$.  With the energy gain-loss formalism of the damped-driven system it is possible to zoom into the transition to chaos using $C = 1/7K$, $C = 1/6K$, and $C = 1/5K$.  We illustrate this transition in Fig. \ref{Fig: Kick Transition} where the system graduates from a sink at $C = 1/7K$, to low period orbits at $C = 1/6K$, to high period seemingly random orbits at $C = 1/5K$, and finally full-scale chaos at $C = 1/K$.

\begin{figure}[htbp]
\vspace{-48pt}
\centering
\stackinset{l}{20mm}{t}{4mm}{\textbf{(a)}}{\includegraphics[width = 0.9\textwidth]{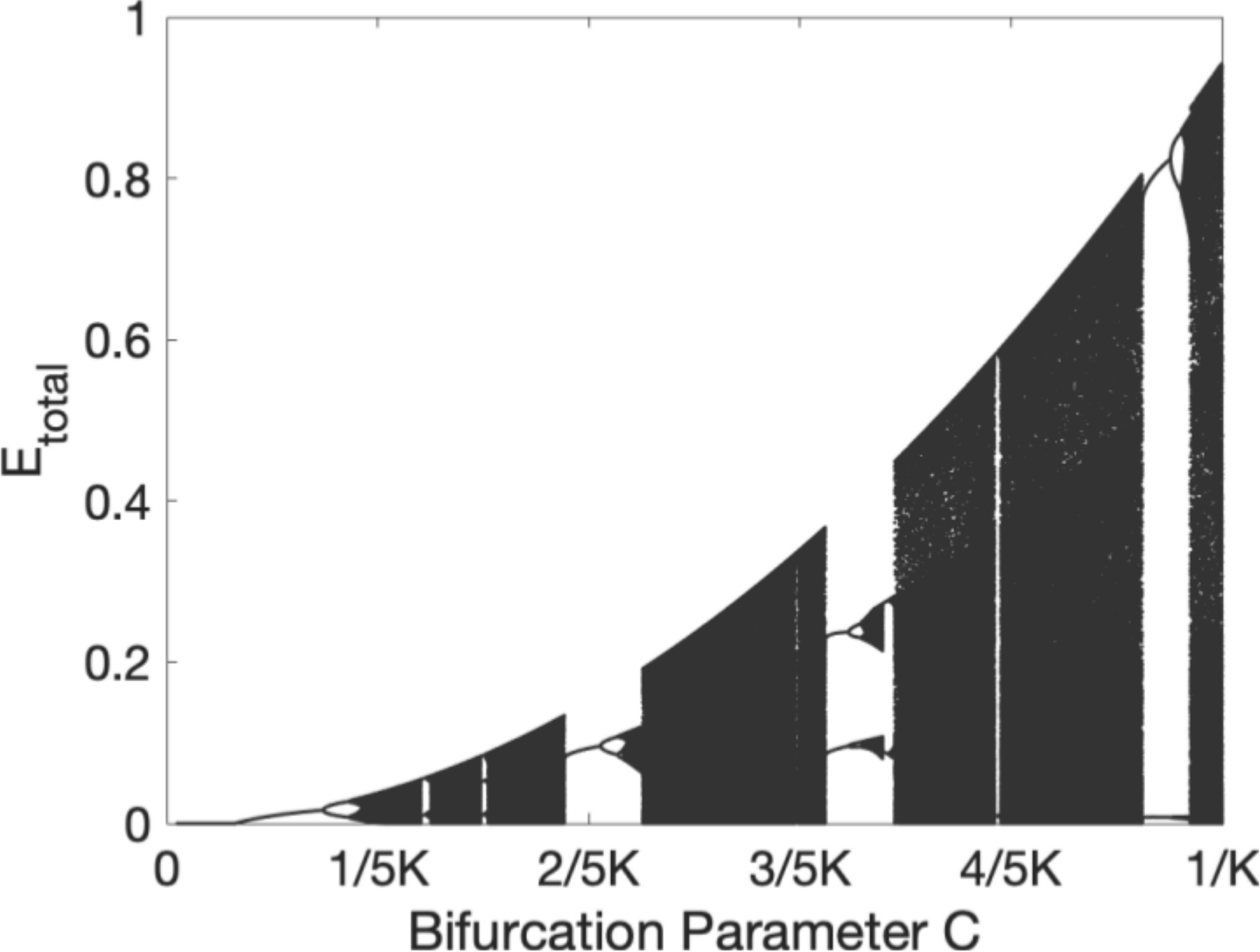}}\qquad
\stackinset{l}{22mm}{t}{4mm}{\textbf{(b)}}{\includegraphics[width = 0.9\textwidth]{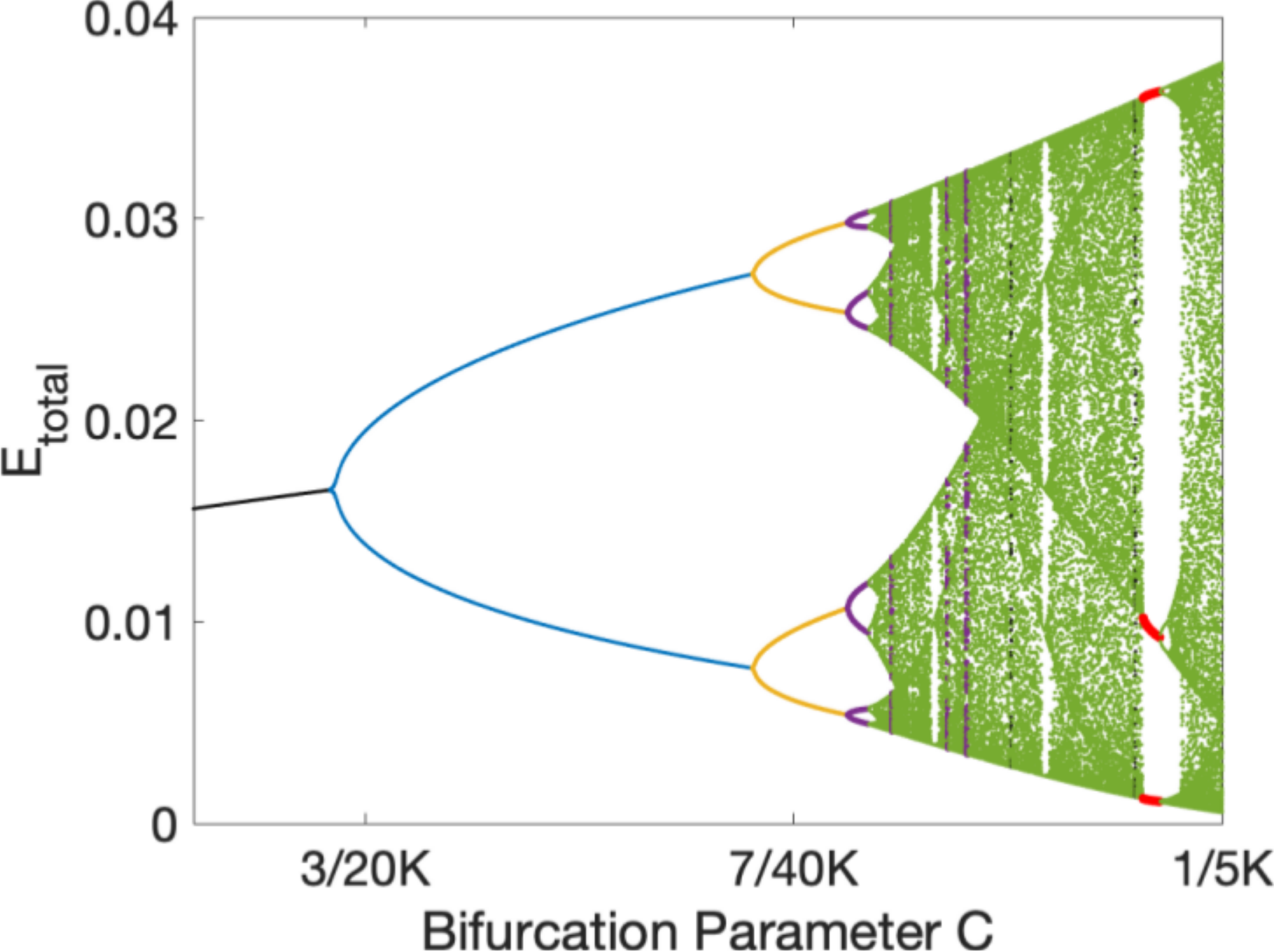}}
\caption{Bifurcation diagram illustrating period doubling bifurcations that give way to orbits of arbitrarily large periods. \textbf{(a)}  The full bifurcation diagram from $C = 0$ to $C = 1/K$.  \textbf{(b)}  A zoomed in bifurcation diagram from $C = 7/50K$ to $C = 1/5K$.  The black curves represent single fixed points, the blue curves represent two cycles, the yellow curves represent four cycles, the purple curves represent eight cycles, and the green curves represent longer periodic orbits.  Finally, the red curves represent three cycles, which in the original 1-D kick map \eqref{Eq: KickModel}, was sufficient to prove chaos \cite{Rahman18}}\label{Fig: Period Doubling}
\end{figure}

A more global view of the bifurcation structure, and hence the route to chaos can be observed in a bifurcation diagram.  It is proved in \cite{Rahman18} that the system \eqref{Eq: KickModel} has a period doubling route to chaos via a flip bifurcation.  The same bifurcation structure is observed in the energy gain-loss formulation of the kick model.  In Fig. \ref{Fig: Period Doubling} we plot the period doubling bifurcation.  As we increase the bifurcation parameter $C$ we observe that the period of the attractors eventually become arbitrarily large, which is an indicator of chaotic dynamics.  While a three cycle is also observed, the proof of Li and Yorke \cite{LiYorke75} only applies to 1-dimensional maps, and therefore in this case is an interesting observation that complements the proof for the 1-D kick map \cite{Rahman18}.

The more geometric bifurcation analysis of the energy gain-loss, compared to the more algebraic analysis of the velocity map in \cite{Rahman18}, reveals windows of seemingly random changes in energy and other windows of energy balance.  This geometric analysis on more complex mechanistic models, which may prohibit the type of analysis conducted in \cite{Rahman18}, will help guide experiments to the approximate location in parameter space of these states.  In doing so we will be able to manipulate the experiments to make stronger connections with systems that rely heavily on the properties of quantum mechanics.  For example, the analysis in this work could guide the creation of an analog for quantum saturable absorbers by using the Energy gain-loss dynamics of mode-locked lasers \cite{LiWaiKutz10} and inferring the necessary physical additions to the annular cavity.

As an example, we show how the energy gain-loss behavior can be inferred from experimental measurements.  Let us consider a droplet at the $n^{\textsuperscript{th}}$ impact.  The kinetic energy of the system just before impact is $E_n^{-} = v_{\text{in}}^2$, where $v_{\text{in}}$ is the measured approach velocity of the droplet relative to the velocity of the table.  Further, the energy right after impact is $E_n^{+} = v_{\text{out}}^2$, where $v_{\text{out}}$ is the measured launch velocity of the droplet relative to the table.  While there is some energy loss during flight, and there is some energy gain during impact, the majority of the energy gain happens during the flight due to gravity and bath acceleration, and a majority of the energy loss happens during impact due to hydrodynamic damping ~\cite{MolBush13a, MolBush13b}.  Therefore, we can partition the dynamics into a gain-regime (the total kinetic energy gained between impacts) and loss-regime (the energy lost due to the impact), which we write as
\begin{subequations}
\begin{align}
E_n^\text{gain} = E_n^{-} - E_{n-1}^{+},\label{Eq: H-K gain}\\
E_n^\text{loss} = E_n^{+} - E_n^{-}. \label{Eq: H-K loss}
\end{align}
\label{Eq: H-K recurrence}
\end{subequations}
Now we have a damped-driven system of the form of Fig. \ref{Fig: Schematic}.  From experimental measurements we can fit the energy gain and energy loss curves, and we can plot an experimentally inferred cobweb to analyze the bifurcations of the damped-driven system.

\section{Comparisons with other damped-driven systems}
\label{Sec: Comparisons}

One main motivation for the present analysis is to facilitate comparisons with other damped-driven systems.  In this section let us briefly compare the energy gain-loss dynamics of \eqref{Eq: Kick Loss-Gain curves} in Sec. \ref{Sec: Damped-driven} to that of mode-locked lasers \cite{LiWaiKutz10} and rotation detonation engines \cite{koch2020modeling}.  
Comparisons can be made through measurements of energy, or proxies thereof.  In optics, for instance, the power spectrum of the electric field provides a proxy measure for the total cavity energy.  Direct measurement of the time-domain evolution of the intensity of the electric field is much more difficult to achieve.  Likewise in rotating detonation engines, the energy expended by the combustion process can either be measured by thermodynamic variables (heat release) or by direct video monitoring of the propagating frame front.  These are indirect measures of the energy expenditures that can be related back directly to the energy expended in the damped-driven system.  

Regardless of the proxy measure, the mode-locked laser and rotating detonation engine both produce observed instabilities consistent with the hydrodynamic quantum analog system.  Specifically, the intersections between the energy gain and loss curves yield cobweb plots that show the transition from period doubling to chaos. In all three systems, this behavior has been observed directly in experiment as well as simulations of the governing PDE models.  

Moreover, repeated non-trivial intersections of the energy gain and loss curves (Fig. \ref{Fig: Kick Transition}), which correspond to the break-up of a pulse for the mode-locked laser \cite{LiWaiKutz10} or having two separate wave fronts for the detonation engines \cite{koch2020modeling}, analogously correspond to a split in the number of preferential velocities for the walking droplet system.  As an example, for mode-locked lasers when the number of intersections between the gain and loss curves go from one to two it indicates a single laser pulse breaking up into two pulses.  For walking droplets, the velocity still varies chaotically, but now instead of being centered about one velocity the distribution will be centered about two velocities.  

Despite the significantly different physical processes of these three systems, optics which is dominated by chromatic dispersion and nonlinear self-focusing, rocket engines which are dominated by combustion and diffusion, and droplets which are dominated by surface tension and collision kinematics with a vibrating fluid bath, the energy balances all produce the same underlying pattern of instability.  Thus damped-driven systems motivate this study of the droplet system in terms of a universal description of a simple gain and loss dynamics and its associated 1-D mapping that we construct explicitly. There are many other damped-driven systems for which this is also the case (bio-locomotion, flutter instabilities, ground-effects for airfoils, etc) and these models will be considered in greater detail in future work.

\section{Conclusion}
\label{Sec: Conclusion}

Damped-driven systems pervade the engineering, physical, and biological sciences.  Despite the diversity of the underlying physics characterizing such systems, they exhibit a common bifurcation structure that is dictated by energy balance considerations.  By separating nonlinear gains and losses, such systems can be represented by compositional mappings.  The wave-particle interaction dynamics of the hydrodynamic quantum analog system is yet another canonical example of a spatio-temporal system whose gain-loss dynamics produce the universal period-doubling cascade to chaos.

Although the governing equations for a droplet interacting with a vibrating fluid bath are quite complex, simplifications of the system \cite{Rahman18} allow us to construct a two-step map which disambiguates losses from gain, or damping from driving.  Our analysis provides explicit approximations for the gain and loss experienced by the droplet system over a short time interval, allowing us to average over spatial effects and construct a compositional map that correctly captures the underlying bifurcations exhibited by the hydrodynamic quantum analog system.  First we treat this superficially through the timeseries data of the model correlated with experimentally validated sample trajectories.  Then we observe the possible intersections between the energy gain and loss curves, which yield cobweb plots that show the transition to chaos.  A non-trivial sink bifurcates into a two-cycle, then orbits of high periodicity, and finally to seemingly random orbits.  Finally, this is verified by bifurcation diagrams of the route to chaos through several period-doubling bifurcations.

The paper also lays the foundational structure upon which the more difficult problem of characterizing the energy gain and loss in a mechanistic hydrodynamic-kinematic model derived directly from experiments can be constructed.  This future work would involve detailed measurements and calculations of the interaction between the droplet and the wave, and also the nonlinear drag due to the intermediate Reynold's number of the air in the presence of the droplet.  Indeed there are many more physical considerations in the full hydrodynamic-kinematic case, but it appears that the underlying dynamics should remain the same.  We see evidence of this in other damped-driven systems \cite{kutz2006mode, koch2020mode} 
and some walking droplet observations \cite{WMHB13, Tambasco2016, BudanurFleury19}.  There are energy gains due to the bath acceleration, and losses due to the hydrodynamic-kinematic interactions.  Competition between the two should lead to the same bifurcation structure presented in this foundational work.

The hydrodynamic quantum analog system is part of a broader class of damped-driven systems, all of which exhibit the same underlying parametric behavior and bifurcation structure.  This includes models from laser physics to rocket engines.  Thus despite the diversity of physics, the compositional mapping acts as a normal form~\cite{nayfeh2011method} encoding for pattern forming instability under very general conditions~\cite{cross1993pattern} for the loss and gain.

\bigskip
\bigskip

\section*{Acknowledgment}

Thanks are due to J.W.M. Bush, D. Harris, and G. Pucci for informative discussions.  A.R. and J.N.K. appreciate the support of the Amath department at UW.  The authors would also like to credit their Amath 499 undergraduate student, L.J. Rhoden for the videos used for Fig. 1.  The work of J.N.K. was also supported in part by the US National Science Foundation (NSF) AI Institute for Dynamical Systems (dynamicsai.org), grant 2112085.  Finally, the authors would like to thank the referees for their detailed suggestions and helpful comments, which were essential in improving the manuscript.

\bigskip
\bigskip

\bibliographystyle{unsrt}
\bibliography{Energy_gain-loss}

\begin{thebibliography}{10}

\bibitem{CPFB05}
Y.~Couder, S.~Protiere, E.~Fort, and A.~Boudaoud.
\newblock Dynamical phenomena: Walking and orbiting droplets.
\newblock {\em Nature}, 437:208, 2005.

\bibitem{CouderFort06}
Y.~Couder and E.~Fort.
\newblock Single-particle diffraction and interference at a macroscopic scale.
\newblock {\em Phys. Rev. Lett.}, 97:154101, 2006.

\bibitem{Bush10}
J.W.M. Bush.
\newblock Quantum mechanics writ large.
\newblock {\em Proc. Nat. Acad. Sci.}, pages 1--2, 2010.

\bibitem{Bush15a}
J.W.M. Bush.
\newblock Pilot-wave hydrodynamics.
\newblock {\em Ann. Rev. Fluid Mech.}, 49:269--292, 2015.

\bibitem{Bush15b}
J.W.M. Bush.
\newblock The new wave of pilot-wave theory.
\newblock {\em Physics Today}, 68(8):47--53, 2015.

\bibitem{BushOza20_ROPP}
J.~W.~M. Bush and A.~U. Oza.
\newblock Hydrodynamic quantum analogs.
\newblock {\em Reports on Progress in Physics}, 84:017001, 2021.

\bibitem{TCB18}
S.~Turton, Miles Couchman, and J.~W.~M. Bush.
\newblock A review of the theoretical modeling of walking droplets: Towards a
  generalized pilot-wave framework.
\newblock {\em Chaos}, 28:096111, 2018.

\bibitem{MolBush13a}
J.~Molacek and J.W.M. Bush.
\newblock Drops bouncing on a vibrating bath.
\newblock {\em J. Fluid Mech.}, 727:582--611, 2013.

\bibitem{MolBush13b}
J.~Molacek and J.W.M. Bush.
\newblock Droplets walking on a vibrating bath: towards a hydrodynamic
  pilot-wave theory.
\newblock {\em J. Fluid Mech.}, 727:612--647, 2013.

\bibitem{HMFCB13}
D.M. Harris, J.~Moukhtar, E.~Fort, Y.~Couder, and J.W.M. Bush.
\newblock Wavelike statistics from pilot-wave dyanmics in a circular corral.
\newblock {\em Phys. Rev. E}, 88:011001, 2013.

\bibitem{OHRB14}
A.~Oza, D.M. Harris, R.R. Rosales, and J.W.M. Bush.
\newblock Pilot-wave dynamics in a rotating frame: on the emergence of orbital
  quantization.
\newblock {\em J. Fluid Mech.}, 744:404--429, 2014.

\bibitem{OWHRB14}
A.~Oza, O.~Wind-Willassen, D.M. Harris, R.R. Rosales, and J.W.M. Bush.
\newblock Pilot-wave dynamics in a rotating frame: Exotic orbits.
\newblock {\em Phys. Fluids}, 26:082101, 2014.

\bibitem{Gilet16}
T.~Gilet.
\newblock Quantumlike statistics of deterministic wave-particle interactions in
  a circular cavity.
\newblock {\em Phys. Rev. E}, 93:042202, 2016.

\bibitem{CSB18}
T.~Cristea-Platon, P.~S\'{a}enz, and J.~W.~M. Bush.
\newblock Walking droplets in a circular corral: Quantization and chaos.
\newblock {\em Chaos}, 28:096116, 2018.

\bibitem{SCB18}
P.~S\'{a}enz, T.~Cristea-Platon, and J.~W.~M. Bush.
\newblock Statistical projection effects in a hydrodynamic pilot-wave system.
\newblock {\em Nature Physics}, 14:315--319, 2018.

\bibitem{EFMC09}
A.~Eddi, E.~Fort, F.~Moisy, and Y.~Couder.
\newblock Unpredictable tunneling of a classical wave-particle association.
\newblock {\em Phys. Rev. Lett.}, 102:240401, 2009.

\bibitem{NMB17}
A.~Nachbin, P.~A. Milewski, and J.~W.~M. Bush.
\newblock Tunneling with a hydrodynamic pilot-wave model.
\newblock {\em Phys. Rev. F}, 2:011001, 2017.

\bibitem{CLE1993}
M.~F. Crommie, C.~P. Lutz, and D.~M. Eigler.
\newblock Confinement of electrons to quantum corrals on a metal surface.
\newblock {\em Science}, 262:218--220, 1993.

\bibitem{QuantumMirage}
H.~C. Manoharan, C.~P. Lutz, and D.~M. Eigler.
\newblock Quantum mirages formed by coherent projection of electronic
  structure.
\newblock {\em Nature}, 403:512--515, 2000.

\bibitem{ORB13}
A.~Oza, R.R. Rosales, and J.W.M. Bush.
\newblock A trajectory equation for walking droplets: hydrodynamic pilot-wave
  theory.
\newblock {\em J. Fluid Mech.}, 737:552--570, 2013.

\bibitem{Gilet14}
T.~Gilet.
\newblock Dynamics and statistics of wave-particle interaction in a confined
  geometry.
\newblock {\em Phys. Rev. E}, 90:052917, 2014.

\bibitem{Rahman18}
Aminur Rahman.
\newblock Standard map-like models for single and multiple walkers in an
  annular cavity.
\newblock {\em Chaos}, 28:096102, 2018.

\bibitem{Neimark}
Ju.I. Neimark.
\newblock On some cases of periodic motions depending on parameters.
\newblock {\em Dokl. Akad. Nauk SSSR}, 129:736--739, 1959.

\bibitem{Sacker}
R.~Sacker.
\newblock On invariant surfaces and bifurcation of periodic solutions of
  ordinary differential equations.
\newblock {\em Report IMM-NYU}, 333:1--62, 1964.

\bibitem{RahmanBlackmore16}
A.~Rahman and D.~Blackmore.
\newblock Neimark-sacker bifurcations and evidence of chaos in a discrete
  dynamical model of walkers.
\newblock {\em Chaos, Solitons \& Fractals}, 91:339--349, 2016.

\bibitem{RJB17}
A~Rahman, Y.~Joshi, and D~Blackmore.
\newblock Sigma map dynamics and bifurcations.
\newblock {\em Regul. Chaotic Dyn.}, 22(6):740--749, 2017.

\bibitem{RahmanBlackmore20}
Aminur Rahman and D.~Blackmore.
\newblock Interesting bifurcations in walking droplet dynamics.
\newblock {\em Commun. Nonlinear Sci. Numer. Simul.}, 90:105348, 2020.

\bibitem{RahmanBlackmoreReview20}
Aminur Rahman and D.~Blackmore.
\newblock Walking droplets through the lens of dynamical systems.
\newblock {\em Mod. Phys. Lett. B}, 34(34):2030009, 2020.

\bibitem{haus2000mode}
Herman~A Haus.
\newblock Mode-locking of lasers.
\newblock {\em IEEE Journal of Selected Topics in Quantum Electronics},
  6(6):1173--1185, 2000.

\bibitem{kutz2006mode}
J~Nathan Kutz.
\newblock Mode-locked soliton lasers.
\newblock {\em SIAM review}, 48(4):629--678, 2006.

\bibitem{ding2009operating}
Edwin Ding and J~Nathan Kutz.
\newblock Operating regimes, split-step modeling, and the haus master
  mode-locking model.
\newblock {\em JOSA B}, 26(12):2290--2300, 2009.

\bibitem{LiWaiKutz10}
Feng Li, P.~K.~A. Wai, and J.~Nathan Kutz.
\newblock Geometrical description of the onset of multi-pulsing in mode-locked
  laser caviities.
\newblock {\em J. Opt. Soc. Am. B}, 27(10):2068 -- 2077, 2010.

\bibitem{bale2009transition}
Brandon~G Bale, Khanh Kieu, J~Nathan Kutz, and Frank Wise.
\newblock Transition dynamics for multi-pulsing in mode-locked lasers.
\newblock {\em Optics express}, 17(25):23137--23146, 2009.

\bibitem{koch2020mode}
James Koch, Mitsuru Kurosaka, Carl Knowlen, and J~Nathan Kutz.
\newblock Mode-locked rotating detonation waves: Experiments and a model
  equation.
\newblock {\em Physical Review E}, 101(1):013106, 2020.

\bibitem{koch2020modeling}
James Koch and J~Nathan Kutz.
\newblock Modeling thermodynamic trends of rotating detonation engines.
\newblock {\em Physics of Fluids}, 32(12):126102, 2020.

\bibitem{koch2021multiscale}
James Koch, Mitsuru Kurosaka, Carl Knowlen, and J~Nathan Kutz.
\newblock Multiscale physics of rotating detonation waves: Autosolitons and
  modulational instabilities.
\newblock {\em Physical Review E}, 104(2):024210, 2021.

\bibitem{agrawal2012fiber}
Govind~P Agrawal.
\newblock {\em Fiber-optic communication systems}, volume 222.
\newblock John Wiley \& Sons, 2012.

\bibitem{kutz1998hamiltonian}
J~Nathan Kutz, Philip Holmes, Stephen~G Evangelides, and James~P Gordon.
\newblock Hamiltonian dynamics of dispersion-managed breathers.
\newblock {\em JOSA B}, 15(1):87--96, 1998.

\bibitem{kutz1999dynamics}
J~Nathan Kutz and Philip Holmes.
\newblock Dynamics and bifurcations of a planar map modeling dispersion managed
  breathers.
\newblock {\em SIAM Journal on Applied Mathematics}, 59(4):1288--1302, 1999.

\bibitem{Tambasco2016}
L.~D. Tambasco, D.~M. Harris, A.~U. Oza, R.~R. Rosales, and J.~W.~M. Bush.
\newblock The onset of chaos in orbital pilot-wave dynamics.
\newblock {\em Chaos}, 26(103107), 2016.

\bibitem{Durey2020c}
M.~Durey.
\newblock Bifurcations and chaos in a {L}orenz-like pilot-wave system.
\newblock {\em Chaos}, 30:103115, 2020.

\bibitem{LogisticMap}
R.~M. May.
\newblock Simple mathematical models with very complicated dynamics.
\newblock {\em Nature}, 26(5560):457, 1976.

\bibitem{Chirikov1979}
B.~V. Chirikov.
\newblock A universal instability of many-dimensional oscillator systems.
\newblock {\em Phys. Rep.}, 52(5):263--379, 1979.

\bibitem{KickedRotator}
B.~V. Chirikov, G.~Casati, F.~M. Izrailev, and J.~Ford.
\newblock {\em Lecture Notes in Physics}, chapter Stochastic behavior of a
  quantum pendulum under a periodic perturbation.
\newblock Springer-Berlin, 1978.

\bibitem{GoodmanMap2008}
R.~H. Goodman.
\newblock Chaotic scattering in solitary wave interactions: A singular
  iterated-map description.
\newblock {\em Chaos}, 18(2):023113, 2008.

\bibitem{GoodmanRahman2015}
R.~H. Goodman, A.~Rahman, M.~Bellanich, and C.~N. Morrison.
\newblock A mechanical analog of the two-bounce resonance of solitary waves:
  Modeling and experiment.
\newblock {\em Chaos}, 25:043109, 2015.

\bibitem{kutz1997mode}
J~Nathan Kutz, Brandon~C Collings, Keren Bergman, Sergio Tsuda, Steven~T
  Cundiff, Wayne~H Knox, Philip Holmes, and Michael Weinstein.
\newblock Mode-locking pulse dynamics in a fiber laser with a saturable bragg
  reflector.
\newblock {\em JOSA B}, 14(10):2681--2690, 1997.

\bibitem{spaulding2002nonlinear}
Kristin~M Spaulding, Darryl~H Yong, Arnold~D Kim, and J~Nathan Kutz.
\newblock Nonlinear dynamics of mode-locking optical fiber ring lasers.
\newblock {\em JOSA B}, 19(5):1045--1054, 2002.

\bibitem{proctor2005passive}
Joshua~L Proctor and J~Nathan Kutz.
\newblock Passive mode-locking by use of waveguide arrays.
\newblock {\em Optics letters}, 30(15):2013--2015, 2005.

\bibitem{proctor2005nonlinear}
Joshua Proctor and J~Nathan Kutz.
\newblock Nonlinear mode-coupling for passive mode-locking: application of
  waveguide arrays, dual-core fibers, and/or fiber arrays.
\newblock {\em Optics express}, 13(22):8933--8950, 2005.

\bibitem{intrachat2003theory}
Karen Intrachat and J~Nathan Kutz.
\newblock Theory and simulation of passive modelocking dynamics using a
  long-period fiber grating.
\newblock {\em IEEE journal of quantum electronics}, 39(12):1572--1578, 2003.

\bibitem{kutz2008passive}
J~Nathan Kutz.
\newblock Passive mode-locking using phase-sensitive amplification.
\newblock {\em Physical Review A}, 78(1):013845, 2008.

\bibitem{FHV15}
B.~Filoux, M.~Hubert, and N.~Vandewalle.
\newblock Strings of droplets propelled by coherent waves.
\newblock {\em Phys. Rev. E}, 92:041004(R), 2015.

\bibitem{WMHB13}
O.~Wind-Willassen, J.~Molacek, D.~Harris, and J.~W.~M. Bush.
\newblock Exotic states of bouncing and walking droplets.
\newblock {\em Phys. Fluids}, 25:082002, 2013.

\bibitem{LiYorke75}
T.Y. Li and J.A. Yorke.
\newblock Period three implies chaos.
\newblock {\em The American Mathematical Monthly}, 82(10):985--992. (DOI:
  10.2307/2318254), 1975.

\bibitem{BudanurFleury19}
N.~B. Budanur and Marc Fleury.
\newblock State space geometry of the chaotic pilot-wave hydrodynamics.
\newblock {\em Chaos}, 29:013122, 2019.

\bibitem{nayfeh2011method}
Ali~H Nayfeh.
\newblock {\em The method of normal forms}.
\newblock John Wiley \& Sons, 2011.

\bibitem{cross1993pattern}
Mark~C Cross and Pierre~C Hohenberg.
\newblock Pattern formation outside of equilibrium.
\newblock {\em Reviews of modern physics}, 65(3):851, 1993.

\end{thebibliography}

\end{document}